\documentclass[10pt,a4paper]{article}
\usepackage[utf8]{inputenc}
\usepackage[T1]{fontenc}
\usepackage{nicefrac}
\usepackage[table]{xcolor}
\usepackage{tikz}
\usetikzlibrary{shapes,arrows,calc, fit, decorations}
\usetikzlibrary{decorations.pathreplacing}
\usetikzlibrary{fit}
\newcommand\addvmargin[1]{
	\node[fit=(current bounding box),inner ysep=#1,inner xsep=0]{};
}

\usepackage{float}
\usepackage{a4wide}
\usepackage{amssymb}

\usepackage{pifont}

\usepackage{mathtools}

\DeclarePairedDelimiter\floor{\lfloor}{\rfloor}

\usepackage[makeroom]{cancel}

\usepackage{array}
\newcolumntype{M}[1]{>{\centering\arraybackslash}m{#1}}
\newcolumntype{N}{@{}m{0pt}@{}}

\usepackage{amsthm}
\usepackage{amsmath}
\usepackage{mathrsfs}

\usepackage{wrapfig}
\usepackage[lofdepth,lotdepth]{subfig}
\usepackage{caption}

\usepackage{stackengine}
\newcommand\xrowht[2][0]{\addstackgap[.5\dimexpr#2\relax]{\vphantom{#1}}}

\usepackage{authblk}

\newtheorem{theorem}{Theorem}
\newtheorem{corollary}[theorem]{Corollary}
\newtheorem{lemma}[theorem]{Lemma}

\newcounter{configSPf}

\newcounter{configSPfp}

\newcounter{configSPn}

\newcounter{configSPqp}

\newcounter{ruleSPf}

\newcounter{ruleSPfp}

\newcounter{ruleSPn}

\newcounter{ruleSPqp}

\title{The chromatic number of 2-edge-colored and signed graphs of bounded maximum degree\footnote{Funding: This work was partially supported by the grant HOSIGRA funded by the French National Research Agency (ANR, Agence Nationale de la Recherche) under the contract number ANR-17-CE40-0022.}}
\author[1]{Christopher Duffy}
\author[2]{Fabien Jacques}
\author[2]{Mickaël Montassier}
\author[2]{Alexandre Pinlou}

\affil[1]{Mathematics and Statistics, University of Saskatchewan, Canada}
\affil[2]{LIRMM, Université de Montpellier, CNRS, Montpellier, France}
\affil[ ]{\footnotesize{\it duffy@math.usask.ca, fabien.jacques@lirmm.fr, mickael.montassier@lirmm.fr, alexandre.pinlou@lirmm.fr}}
\begin{document}

\maketitle

\begin{abstract}
A 2-edge-colored graph or a signed graph is a simple graph with two types of edges. A homomorphism from a 2-edge-colored graph $G$ to a 2-edge-colored graph $H$ is a mapping $\varphi: V(G) \rightarrow V(H)$ that maps every edge in $G$ to an edge of the same type in $H$. Switching a vertex $v$ of a 2-edge-colored or signed graph corresponds to changing the type of each edge incident to $v$. There is a homomorphism from the signed graph $G$ to the signed graph $H$ if after switching some subset of the vertices of $G$ there is a 2-edge-colored homomorphism from $G$ to $H$. 

The chromatic number of a 2-edge-colored (resp. signed) graph $G$ is the order of a smallest 2-edge-colored (resp. signed) graph $H$ such that there is a homomorphism from $G$ to $H$. The chromatic number of a class of graph is the maximum of the chromatic numbers of the graphs in the class.

We study the chromatic numbers of 2-edge-colored and signed graphs (connected and not necessarily connected) of a given bounded maximum degree. More precisely, we provide exact bounds for graphs of maximum degree 2. We then propose specific lower and upper bounds for graphs of maximum degree 3, 4, and 5. We finally propose general bounds for graphs of maximum degree $k$, for every $k$.

\end{abstract}


\section{Introduction}

\subsection{Signed and 2-edge-colored graphs}

A \textit{2-edge-colored graph} or a \textit{signed graph} $G = (V, E, s)$ is a simple graph $(V, E)$ with two kinds of edges: positive and negative edges. We do not allow parallel edges nor loops. The signature $s: E(G) \rightarrow \{-1, +1\}$ assigns to each edge its sign. For the concepts discussed in this article, 2-edge-colored graphs and signed graphs only differ on the notion of homomorphism.
Note that 2-edge-colored graphs are sometimes referred to as \textit{signified} graphs by some authors.

A \textit{positive neighbor} (resp. \textit{negative neighbor}) of a vertex $v$ is a vertex that is connected to $v$ with a positive (resp. negative) edge.

\textit{Switching} a vertex $v$ of a $2$-edge-colored or signed graph corresponds to reversing the signs of all the edges that are incident to $v$. 

Two 2-edge-colored or signed graphs $G$ and $G'$ are \textit{switching equivalent} if it is possible to turn $G$ into $G'$ after any number of switches.

Given a 2-edge-colored or signed graph $G = (V, E, s)$, the \textit{underlying graph} of $G$ is the simple graph $(V, E)$.

A cycle of a $2$-edge-colored or signed graph is said to be \textit{balanced} (resp. \textit{unbalanced}) if it has an even (resp. odd) number of negative edges. 
The notion of balanced cycles allows us to define switching equivalence as follows.

\begin{theorem}[Zaslavsky \cite{Z}]
\label{thm:Zaslavsky}
Two $2$-edge-colored or signed graphs are switching equivalent if and only if they have the same underlying graph and the same set of balanced cycles.
\end{theorem}

\subsection{Homomorphisms}
Given two $2$-edge-colored graphs $G$ and $H$, the mapping $\varphi : V(G)\rightarrow V(H)$ is a \textit{homomorphism} if $\varphi$ maps every edge of $G$ to an edge of $H$ with the same sign. This can be seen as coloring the vertices of $G$ by using the vertices of $H$ as colors. The target graph $H$ gives us the rules that this coloring must follow. If vertices $1$ and $2$ of $H$ are adjacent with a positive (resp. negative) edge, then every pair of adjacent vertices in $G$ colored with $1$ and $2$ must be adjacent with a positive (resp. negative) edge. 

If $G$ admits a homomorphism to $H$, we say that $G$ is \textit{$H$-colorable} or that $H$ can color $G$. If $G$ admits a homomorphism to a graph on $n$ vertices, we say that $G$ is \textit{$n$-colorable}.

The \textit{chromatic number} $\chi_2(G)$ of a $2$-edge-colored graph $G$ is the order (the number of vertices) of a smallest $2$-edge-colored graph $H$ such that $G$ is $H$-colorable. 
The chromatic number $\chi_2(\mathcal{C})$ of a class of 2-edge-colored graphs $\mathcal{C}$ is the maximum of the chromatic numbers of the graphs in the class.

A \textit{2-edge-colored clique} is a 2-edge-colored graph that has the same order and chromatic number.

\begin{lemma}[\cite{B17}]
\label{lem:2ecC}
A 2-edge-colored graph is a 2-edge-colored clique if and only if each pair of non-adjacent vertices is connected by a path of length 2 made of one positive and one negative edge.
\end{lemma}

Given two signed graphs $G$ and $H$, the mapping $\varphi : V(G)\rightarrow V(H)$ is a \textit{homomorphism} if there is a homomorphism from $G$ to $H$ after switching some subset of the vertices of $G$ and/or switching some subset of the vertices of $H$. However, switching in $H$ is unnecessary (as explained in Section 3.3 of \cite{HomSG}).

The \textit{chromatic number} $\chi_s(G)$ of a signed graph $G$ is the order of a smallest signed graph $H$ such that $G$ admits a homomorphism to $H$. 
The chromatic number $\chi_s(\mathcal{C})$ of a class of signed graphs $\mathcal{C}$ is the maximum of the chromatic numbers of the graphs in the class. \newline

A \textit{signed clique} is a signed graph that has the same order and chromatic number.
\begin{lemma}[\cite{HomSG}]
\label{lem:sC}
A signed graph is a signed clique if and only if every pair of non-adjacent vertices is part of an unbalanced cycle of length 4.
\end{lemma}

A class of graphs is \textit{colorable} if there exists a target graph that can color every graph in the class.
A class of graphs is \textit{complete} if for every two graphs $G_1$ and $G_2$ in the class, there is a graph $G^*$ in the class such that $G_1$ and $G_2$ are subgraphs of $G^*$.

A class $\mathcal{C}$ of 2-edge-colored (resp. signed) graphs is \textit{optimally colorable} if there exists a target 2-edge-colored (resp. signed) graph $T$ on $\chi_2(\mathcal{C})$ (resp.  $\chi_s(\mathcal{C})$ ) vertices such that every graph in $\mathcal{C}$ is $T$-colorable. 

\begin{lemma}
\label{lem:opt_col}
Every class $\mathcal{C}$ of graphs which is colorable and complete is optimally colorable.
\end{lemma}

\begin{proof}
Following the proof in \cite{ChiOG}: Suppose that $\mathcal{C}$ is colorable and complete but not optimally colorable. There exists a finite set $S$ of graphs in $\mathcal{C}$ which cannot be colored with a single target graph on $\chi(\mathcal{C})$ vertices (such a set can be finite since there exists a finite number of target graphs having at most $\chi(\mathcal{C})$ colors). Since $\mathcal{C}$ is complete, there exists a graph $G$ in $\mathcal{C}$ that contains every graph in $S$ as subgraphs. Graph $G$ admits a homomorphism to a target graph $T$ on $\chi(\mathcal{C})$ vertices. Therefore, every graph in $S$ can be colored with $T$, a contradiction.
\end{proof}

$2$-edge-colored graphs are, in some sense, similar to oriented graphs since a pair of vertices can be adjacent in two different ways in both kinds of graphs: with a positive or a negative edge in the case of $2$-edge-colored graphs, with a toward or a backward arc in the oriented case. 

The notion of homomorphism of oriented graphs has been introduced by Courcelle \cite{Cou94} in 1994 and has been widely studied since then. Due to the similarity above-mentioned, we try to adapt techniques used to study the homomorphisms of oriented graphs of bounded degree to 2-edge-colored graphs of bounded degree.

\subsection{Target Graphs}
A $2$-edge-colored graph $(V, E, s)$ is said to be \textit{antiautomorphic} if it is isomorphic to $(V, E, -s)$.

A $2$-edge-colored graph $G = (V, E, s)$ is said to be \textit{$K_n$-transitive} if for every pair of cliques $\{u_1, u_2, \ldots , u_n\}$ and $\{v_1, v_2, \ldots , v_n\}$ in $G$ such that $s(u_i u_j) = s(v_i v_j)$ for all $i \neq j$, there exists an automorphism that maps $u_i$ to $v_i$ for all $i$. 
For $n = 1$, $2$, or $3$, we say that the graph is \textit{vertex-transitive}, \textit{edge-transitive}, or \textit{triangle-transitive}, respectively. 

A 2-edge-colored graph $G$ has \emph{Property $P_{k, n}$} if for every sequence of $k$ distinct vertices $(v_1, v_2, \dots, v_k)$ that induces a clique in $G$ and for every sign vector $(\alpha_1, \alpha_2, ..., \alpha_k) \in \{-1, +1\}^k$ there exist at least $n$ distinct vertices $\{u_1, u_2, ..., u_n\}$ such that $s(v_i u_j) = \alpha_i$ for $1 \leq i \leq k$ and $1 \leq j \leq n$.

Let $q$ be a prime power with $q \equiv 1 \mod 4$. Let $\mathbb{F}_q$ be the finite field of order $q$.

The \textit{2-edge-colored Paley graph} $SP_q$ has vertex set $V(SP_q) = \mathbb{F}_q$. Two vertices $u$ and $v \in V(SP_q)$, $u \neq v$, are connected with a positive edge if $u - v$ is a square in $\mathbb{F}_q$ and with a negative edge otherwise.
This definition is consistent since $q \equiv 1 \mod 4$ so $-1$ is always a square in $\mathbb{F}_q$ and if $u-v$ is a square then $v-u$ is also a square.

\begin{lemma}[\cite{OPS16}]
\label{lem:PSP}
Graph $SP_q$ is vertex-transitive, edge-transitive, antiautomorphic and has Properties $P_{1, \frac{q-1}{2}}$ and $P_{2, \frac{q-5}{4}}$.
\end{lemma}

Given a 2-edge-colored graph $G$ with signature $s_G$, we create the \textit{antitwinned graph} of $G$ denoted by $\rho(G)$ as follows. \newline

Let $G^{+1}$, $G^{-1}$ be two copies of $G$. The vertex corresponding to $v \in V(G)$ in $G^{i}$ is denoted by~$v^i$. The vertex set, edge set and signature of $\rho(G)$ is defined as follows:

\begin{itemize}
\item $V(\rho(G)) = V(G^{+1}) \cup V(G^{-1})$
\item $E(\rho(G)) = \{ u^i v^j : uv \in E(G), \ i, j \in \{-1, +1\} \}$
\item $s_{\rho(G)}(u^i v^j) = i \times j \times s_G(uv)$
\end{itemize}

By construction, for every vertex $v$ of $G$, $v^{-1}$ and $v^{+1}$ are \textit{antitwins}, the positive neighbors of $v^{-1}$ are the negative neighbors of $v^{+1}$ and vice versa. A 2-edge-colored graph is \textit{antitwinned} if every vertex has a unique antitwin.

\begin{lemma}[\cite{HomEG}]
\label{lem:BG}
Let $G$ and $H$ be 2-edge-colored graphs. The two following propositions are equivalent:
\begin{itemize}
\item The graph $G$ admits a homomorphism to $\rho(H)$.
\item The graph $G$, seen as a signed graph, admits a homomorphism to $H$.
\end{itemize}
\end{lemma}

In other words, if a 2-edge-colored graph admits a homomorphism to an antitwinned target graph on $n$ vertices,  then the same graph as seen as a signed graph also admits a homomorphism to a target graph on $\frac{n}{2}$ vertices.
The family $\rho(SP_q)$ are interesting target graphs (especially for bounding the chromatic number of signed graphs since they are antitwinned graphs).

\begin{lemma}[\cite{OPS16}]
\label{lem:PrhoSP}
The graph $\rho(SP_q)$ is vertex-transitive, antiautomorphic and has Properties $P_{1, q-1}$, $P_{2, \frac{q-3}{2}}$ and $P_{3, max(\frac{q-9}{4}, 0)}$.
\end{lemma}

One last family of interesting target graphs are the Tromp-Paley graphs $TR(SP_q)$.
Let $SP_q^+$ be $SP_q$ with an additional vertex that is connected to every other vertex with a positive edge.
The Tromp-Paley graph $TR(SP_q)$ is equal to $\rho(SP_q^+)$.

This construction improves the properties of $\rho(SP_q)$ at the cost of having two more vertices. Since Tromp-Paley graphs are antitwinned, they are interesting for bounding the chromatic number of signed graphs.

\begin{lemma}[\cite{OPS16}]
\label{lem:PTRSP}
The graph $TR(SP_q)$ is vertex-transitive, edge-transitive, antiautomorphic and has properties $P_{1, q}$, $P_{2, \frac{q-1}{2}}$ and $P_{3, \frac{q-5}{4}}$.
\end{lemma}

\section{Results}

In the sequel, $\mathcal{D}_k$ (resp. $\mathcal{D}_k^c$) denotes the class of $2$-edge-colored or signed graphs (resp. connected $2$-edge-colored or signed graphs) with maximum degree $k$, graphs in which a vertex cannot be adjacent to more than $k$ other vertices.

Tables~\ref{table:results2} and \ref{table:resultsS} summarize results on the chromatic number of the classes of (connected) 2-edge-colored and signed graphs of bounded degree. Grey cells contain our results presented in this paper, while white cells contain already known results.
  
\newcommand{\greycell}{\cellcolor{black!20}}

\begin{table}[!h]
\centering
	\begin{tabular}{|c||c|c|c|c|}
		\hline \xrowht[()]{9pt}
		& $\chi_{2}(\mathcal{D}_k)$ & $\chi_{2}(\mathcal{D}_k^c)$ \\
		\hline 
		\hline \xrowht[()]{9pt}
		$k=1 $ & \greycell$\chi_2(\mathcal{D}_1)=3$ & \greycell$\chi_2(\mathcal{D}_1^c)=2$ \\
		\hline \xrowht[()]{9pt}
		$k=2 $ & \greycell$\chi_2(\mathcal{D}_2) = 6$ & \greycell$\chi_2(\mathcal{D}_2^c) = 5$ \\
		\hline \xrowht[()]{9pt}
		$k=3 $ & \greycell$8 \leq \chi_2(\mathcal{D}_3) \leq 11$ & \greycell$8 \leq \chi_2(\mathcal{D}_3^c) \leq 10$ \\
		\hline \xrowht[()]{9pt}
		$k=4 $ & \multicolumn{2}{c|}{\greycell$12 \leq \chi_2(\mathcal{D}_4^c) \leq \chi_2(\mathcal{D}_4) \leq 30	$} \\
		\hline \xrowht[()]{9pt}
		$k=5 $ & \multicolumn{2}{c|}{\greycell$16 \leq \chi_2(\mathcal{D}_5^c) \leq \chi_2(\mathcal{D}_5) \leq 110$} \\
		\hline \xrowht[()]{9pt}
		$6 \leq k \leq 10$ & \greycell$4(k-1) \leq \chi_{2}(\mathcal{D}_k) \leq k^2\cdot2^{k+1}$ &  $4(k-1) \leq \chi_{2}(\mathcal{D}_k^c) \leq (k-1)^2\cdot2^{k}+2$ \cite{MIXED} \\
		\hline \xrowht[()]{9pt}
		$11 \leq k $ & \greycell $2^{\frac{k}{2}} \leq \chi_{2}(\mathcal{D}_k) \leq k^2\cdot2^{k+1}$ &  $2^{\frac{k}{2}} \leq \chi_{2}(\mathcal{D}_k^c) \leq (k-1)^2\cdot2^{k}+2$ \cite{MIXED} \\
		\hline
	\end{tabular}
	\caption{Results on the chromatic number of the classes of (connected) 2-edge-colored graphs of bounded degree.}\label{table:results2}
\end{table}

\begin{table}[!h]
\centering
	\begin{tabular}{|c||c|c|c|c|}
		\hline \xrowht[()]{9pt}
		& $\chi_{s}(\mathcal{D}_k)$ & $\chi_{s}(\mathcal{D}_k^c)$ \\
		\hline 
		\hline \xrowht[()]{9pt}
		$k=1 $ &  \multicolumn{2}{c|}{\greycell$\chi_s(\mathcal{D}_1) = \chi_s(\mathcal{D}_1^c) = 2$} \\
		\hline \xrowht[()]{9pt}
		$k=2 $ & \multicolumn{2}{c|}{\greycell$\chi_s(\mathcal{D}_2) = \chi_s(\mathcal{D}_2^c) = 4$}  \\
		\hline \xrowht[()]{9pt}
		$k=3 $ & $6 \leq \chi_s(\mathcal{D}_3) \leq 7$ \cite{BPS} & $\chi_s(\mathcal{D}_3^c) = 6$ \cite{BPS} \\
		\hline \xrowht[()]{9pt}
		$k=4 $ & \multicolumn{2}{c|}{\greycell$10 \leq \chi_s(\mathcal{D}_4^c)\leq \chi_s(\mathcal{D}_4) \leq 16$} \\
		\hline \xrowht[()]{9pt}
		$k=5 $ & \multicolumn{2}{c|}{\greycell$12 \leq \chi_s(\mathcal{D}_5^c)\leq \chi_s(\mathcal{D}_5) \leq 56$} \\
		\hline \xrowht[()]{9pt}
		$6 \leq k \leq 8$ & \greycell $2(k+1) \leq \chi_{s}(\mathcal{D}_k) \leq k^2\cdot2^{k+1}$ & $2(k+1) \leq \chi_s(\mathcal{D}_k^c) \leq (k-1)^2\cdot2^{k}+2$ \cite{MIXED} \\
		\hline \xrowht[()]{9pt}
		$9 \leq k $ & \greycell $2^{\frac{k}{2}-1} \leq \chi_{s}(\mathcal{D}_k) \leq k^2\cdot2^{k+1}$ & $2^{\frac{k}{2}-1} \leq \chi_s(\mathcal{D}_k^c)\leq (k-1)^2\cdot2^{k}+2$ \cite{MIXED} \\
		\hline
	\end{tabular}
	\caption{Results on the chromatic number of the classes of (connected) signed graphs of bounded degree.}\label{table:resultsS}
\end{table}

An edge of a $2$-edge-colored graph has chromatic number 2 and thus $\chi_2(\mathcal{D}_1^c)=2$; however, a $2$-edge-colored graph with two non-adjacent edges, one positive and one negative, has chromatic number 3 (the target graph needs a positive and a negative edge, hence at least three vertices) and thus $\chi_2(\mathcal{D}_1)=3$. We therefore have a difference between the chromatic numbers of connected and non-connected 2-edge-colored graphs with maximum degree 1. This difference does not exist for signed graphs since a negative edge can be changed into a positive one after a switch.
This difference  (and lack thereof for signed graphs) appears also in graphs with maximum degree 2 as explained in the next subsection. We do not know yet if this is also the case for graphs with maximum degree at least 3.

The rest of this article is devoted to the proofs of the results presented in Table~\ref{table:results2} and \ref{table:resultsS}.

\section{Lower bounds}

We begin with the following theorems that gives us lower bounds for the chromatic numbers of 2-edge-colored and signed graphs.

\begin{theorem}
\label{thm:2ecc}
For every $k \geq 3$ there is a $k$-regular 2-edge-colored clique on $4 \cdot (k-1)$ vertices.
\end{theorem}
\begin{proof}
Let $G$ be the 2-edge-colored graph with vertex set $V(G) = \{0, 1, ..., 4 \cdot (k-1) - 1\}$. In this proof, every number is considered modulo $4 \cdot (k-1) $.
For all $u \in V(G)$:
\begin{itemize}
\item If $u$ is even, $u$ is positively adjacent to $u + 2 \cdot (k-1) $ and $u + 2i + 1 $ for $0 \leq i \leq k-3$ and negatively adjacent to $u - 1 $.
\item If $u$ is odd, $u$ is positively adjacent to $u- 2i - 1 $ for $0 \leq i \leq k-3$ and negatively adjacent to $u + 1 $ and $u + 2(k-1) $.
\end{itemize}

Graph $G$ is $k$-regular. We now show that every pair of vertices is either adjacent or is connected by a path of length 2 made of one positive and one negative edge in order to conclude with Lemma~\ref{lem:2ecC}. It suffices to show that this is the case for each pair of vertices containing $0$ or $1$ (since adding $2 $ to every vertex yields an automorphism).

Vertex $0$ is adjacent to $2(k-1)$, $2i+1$ for $0 \leq i \leq k-3$, $4(k-1)-1$. The following paths are made of one positive and one negative edge: $(0, 2(k-1), 2(k-1)-1)$, $(0, 2i+1, 2i+2)$, $(0, 2i+1, 2(k-1) + 2i +1)$, $(0, 4(k-1)-1, 4(k-1)-2-2i)$ for $0 \leq i \leq k-3$. We have covered all pairs $(0, v)$ with $v \in \{2i+1, 2i+2, 2(k-1)-1, 2(k-1), 2(k-1)+2i+1, 4(k-1)-2-2i, 4(k-1)-1\} = V(G) \setminus \{0\}$.

Vertex $1$ is adjacent to $4(k-1)-2i$ for $0 \leq i \leq k-3$, $2$ and $2(k-1)+1$. The following paths are made of one positive and one negative edge: $(1, 4(k-1)-2i, 4(k-1)-2i-1)$, $(1, 2, 2i+3)$, $(1, 2(k-1)+1, 2(k-1)-2i)$ and $(1, 2(k-1)+1, 2(k-1)+2)$ for $0 \leq i \leq k-3$. We have covered all pairs $(1, v)$ with $v \in \{2, 2i+3, 2(k-1)-2i, 2(k-1)+1, 2(k-1)+2, 4(k-1)-2i-1, 4(k-1)-2i\}= V(G) \setminus \{1\}$.
\end{proof}

\begin{theorem}
\label{thm:sc}
For every $k \geq 4$ there is a $k$-regular signed clique on $2 \cdot (k+1)$ vertices.
\end{theorem}
\begin{proof}
Let $G$ be the signed graph with vertex set $V(G) = \{0, 1, ..., 2 \cdot (k+1)-1\}$. In this proof, every number is considered modulo $2 \cdot (k+1)$.
For all $u \in V(G)$:
\begin{itemize}
\item If $u$ is even, $u$ is positively adjacent to $u+1 $ and $u+4+2i $ for $0 \leq i \leq k-3$ and negatively adjacent to $u -1 $.
\item If $u$ is odd, $u$ is negatively adjacent to $u+1 $ and $u+4+2i $ for $0 \leq i \leq k-3$ and positively adjacent to $u -1 $.
\end{itemize}

Graph $G$ is $k$-regular. We now show that every pair of vertices is part of an unbalanced cycle of length $4$ in order to conclude with Lemma~\ref{lem:sC}. It suffices to show that this is the case for each pair of vertices containing $0$ (since adding $2 $ to every vertex yields an automorphism and adding $1 $ to every vertex yields an antiautomorphism).

Cycles $(0, 1, 2, 2(k+1)-4)$, $(0, 2(k+1)-1, 2(k+1)-2, 4)$, $0, 4, 3, 2(k+1)-1)$, $(0, 4+2i, 5+2i, 1)$ and $(0, 4+2i, 3+2i, 2(k+1)-1)$ for $0 \leq i \leq k-3$ are unbalanced. This covers all pairs $(0, v)$ with $v \in \{1, 2, 3, 4, 3+2i, 4+2i, 5+2i, 2(k+1)-2, 2(k+1)-1\} = V(G) \setminus \{0\}$.
\end{proof}

\section{Graphs with maximum degree 2}
This section is devoted to 2-edge-colored and signed graphs with maximum degree $2$. We prove that $\chi_2(\mathcal{D}_2^c) = 5$, $\chi_2(\mathcal{D}_2) = 6$ and $\chi_s(\mathcal{D}_2) = \chi_s(\mathcal{D}_2^c) = 4$.
\subsection{Connected $2$-edge-colored graphs with maximum degree $2$}
\label{sec:2c}

In this subsection, we consider the case of connected $2$-edge-colored graphs with maximum degree~$2$ and we prove that their chromatic number is exactly $5$. We obtain this result by showing that every graph $G\in\mathcal{D}_2^c$ admits a homomophism to one of the two graphs of Figure~\ref{fig:2c_targets}.

\begin{theorem}[$\chi_2(\mathcal{D}_2^c) = 5$]
\label{thm:2c}
	The class of connected $2$-edge-colored graphs with maximum degree $2$ has chromatic number $5$ and is not optimally colorable. 
\end{theorem}

\begin{proof}
The class of connected graphs with maximum degree 2 is the set of all paths and cycles. The cycle of length 6 from Figure~\ref{fig:2_counter-example} has chromatic number 5. We start by showing that it is not possible to color it with four colors.\newline

Vertices $v_1$, $v_2$ and $v_3$ belong to a path of length 2 with one negative and one positive edge. We therefore need 3 distinct colors for these vertices and without loss of generality we color $v_1$, $v_2$ and $v_3$ with $1$, $2$ and $3$ respectively. Using the same argument, $v_4$ cannot receive colors $2$ or $3$.

Suppose that we color $v_4$ in $1$. Vertex $v_5$ cannot be colored in $1$, $2$ or $3$ so we color it in $4$. We would need a new color to color $v_6$.

Suppose that we color $v_4$ in $4$. Vertex $v_5$ cannot be colored in $3$ or $4$. If we color $v_5$ in $1$ it will not be possible to color $v_6$. If we color $v_5$ in $2$ we would need a new color to color $v_6$. \newline

Therefore, it is not possible to color this graph with $4$ colors. A 5-coloring exists (we color the vertices with $1$, $2$, $3$, $4$, $5$ and $3$ in order) so the chromatic number of this 2-edge-colored graph is 5 and the class of connected $2$-edge-colored graphs with maximum degree 2 has chromatic number at least 5. \newline

We now show that any 2-edge-colored graph with maximum degree two admits a homomorphism to $SP_5$ (see Figure~\ref{sub:SP5}), the signed Paley graph on 5 vertices, or $SB$, the signed butterfly (see Figure~\ref{sub:SB}).

\begin{figure}[H]
  \begin{center}
    \subfloat[The 2-edge-colored graph $SP_5$.]{
    	\scalebox{0.7}
		{
			\begin{tikzpicture}[thick]
			\def \radius {3cm}
			\def \margin {8} 
			\tikzstyle{vertex}=[circle,minimum width=2.5em]
			
			\node[draw, vertex] (0) at (-72*0+90: 2) {$0$};
    		\node[draw, vertex] (1) at (-72*1+90: 2) {$3$};
    		\node[draw, vertex] (2) at (-72*2+90: 2) {$4$};
    		\node[draw, vertex] (3) at (-72*3+90: 2) {$2$};
    		\node[draw, vertex] (4) at (-72*4+90: 2) {$1$};
    		
    		\draw (0) -- (1);
    		\draw (1) -- (2);
    		\draw (2) -- (3);
    		\draw (3) -- (4);
    		\draw (4) -- (0);
    		
    		\draw[dashed] (0) -- (2);
    		\draw[dashed] (1) -- (3);
    		\draw[dashed] (2) -- (4);
    		\draw[dashed] (3) -- (0);
    		\draw[dashed] (4) -- (1);
			\end{tikzpicture}
			
		}
      \label{sub:SP5}
    } \hspace*{3cm}
    \subfloat[The 2-edge-colored graph $SB$.]{
    \scalebox{0.8}
		{
			\begin{tikzpicture}[thick]
			\def \radius {3cm}
			\def \margin {8} 
			\tikzstyle{vertex}=[circle,minimum width=1.5em]
			
			\node[draw, vertex] (0) at (0: 0) {};
    		\node[draw, vertex] (1) at (30: 2) {};
    		\node[draw, vertex] (2) at (-30: 2) {};
    		\node[draw, vertex] (3) at (30+180: 2) {};
    		\node[draw, vertex] (4) at (-30+180: 2) {};
    		
    		\draw (0) -- (1);
    		\draw (1) -- (2);
    		\draw (2) -- (0);
    		
    		\draw[dashed] (0) -- (3);
    		\draw[dashed] (3) -- (4);
    		\draw[dashed] (4) -- (0);
			\end{tikzpicture}
			
		}
      \label{sub:SB}
    }
    \caption{Every connected $2$-edge-colored graph with maximum degree 2 can be colored with at least one of these two graphs.} 
    \label{fig:2c_targets}  
  \end{center}
\end{figure}
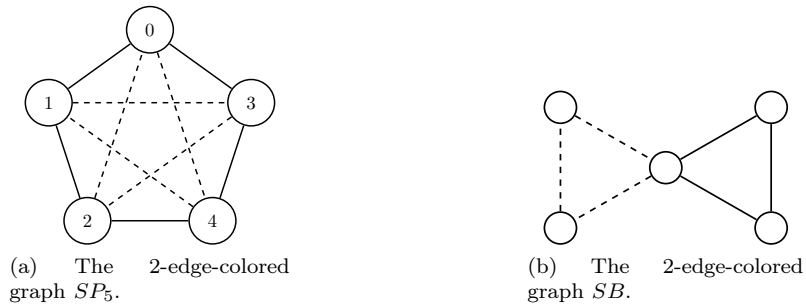

Notice that any 2-edge-colored path can be colored with the graph from Figure~\ref{fig:2c_paths_target} because every vertex in this graph has at least one positive and at least one negative neighbor. This graph is a subgraph of $SP_5$, thus every path maps to $SP_5$. In the following we refer to vertices with even or odd labels as even or odd vertices. Note that in this subgraph, odd (resp. even) vertices are only connected to even (resp. odd) vertices. Also note that every odd (resp. even) vertex of this subgraph is linked with a positive (resp. negative) edge to 0 in $SP_5$.

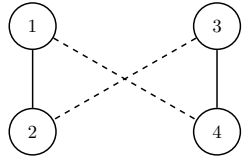
\begin{figure}[H]
	\centering
	\scalebox{0.7}
	{
		\begin{tikzpicture}[thick]
		\def \radius {3cm}
		\def \margin {8} 
		\tikzstyle{vertex}=[circle,minimum width=2.5em]
		
		\node[draw, vertex] (3) at (-60*1+90: 2) {$3$};
		\node[draw, vertex] (2) at (-60*2+90: 2) {$4$};
		\node[draw, vertex] (4) at (-60*4+90: 2) {$2$};
		\node[draw, vertex] (1) at (-60*5+90: 2) {$1$};
		
		\draw[dashed] (1) -- (2);
		\draw[dashed] (3) -- (4);
		
		\draw (1) -- (4);
		\draw (2) -- (3);
		\end{tikzpicture}
		
	}
	\caption{Target graph that can color any path.}
	\label{fig:2c_paths_target}
\end{figure}

Let $G = (V, E, s)$ be a 2-edge-colored cycle with $V(G) = \{v_0, v_1, ..., v_{n-1}\}$ and $E(G) = \{v_i v_j | i-j \equiv 1 \mod n\}$. We now create a homomorphism $\varphi$ from $G$ to $SP_5$ or $SB$. \newline

Suppose that $n$ is even:

Suppose there is a vertex which is incident to two positive edges. Without loss of generality, let $v_0$ be this vertex. We create $\varphi : G \rightarrow SP_5$ as follows. Color $\varphi(v_0)$ is equal to 0. We then color the path $\{v_1, v_2, ..., v_{n-1}\}$ with the subgraph from Figure~\ref{fig:2c_paths_target}. Since $s(v_0 v_1) = +1$, $\varphi(v_1)$ has to be an odd color. Since every odd (resp. even) vertex of the subgraph is only adjacent to even (resp. odd) vertices, we alternate between odd and even colors along the path $\{v_2, v_3, ..., v_{n-1}\}$. Hence $v_{n-1}$ is colored in an odd color and completes the homomorphism since $s(v_{n-1} v_0) = +1$.

Similarly, if there is a vertex which is incident to two negative edges, we can also create a homomorphism $\varphi : G \rightarrow SP_5$. 

We can now assume that the cycle alternates between positive and negative edges. Without loss of generality let $s(v_0 v_1) = -1$. We create $\varphi : G \rightarrow SP_5$ as follows:

$$
\varphi(v_i) = \left\{
    \begin{array}{ll}
        0 & \mbox{if } i = 0, \\
        4 & \mbox{if } i = 1, \\
        2 & \mbox{if } i \equiv 2 \mod 4, \\
        3 & \mbox{if } i \equiv 3 \mod 4, \\
        4 & \mbox{if } i \equiv 0 \mod 4 \mbox{ and } i \neq 0, \\
        1 & \mbox{if } i \equiv 1 \mod 4 \mbox{ and } i \neq 1. \\
    \end{array}
\right.
$$

The color of $v_{n-1}$ will thus be an odd color and complete the homomorphism since $s(v_{n-1} v_0) = +1$. \newline

Suppose that $n$ is odd:

Suppose there is a vertex which is incident to one positive and one negative edge. Without loss of generality let $s(v_{n-1} v_0) = -1$ and $s(v_0 v_1) = +1$. We create $\varphi : G \rightarrow SP_5$ as follows. Color $\varphi(v_0)$ is equal to 0. We then color the path $\{v_1, v_2, ..., v_{n-1}\}$ with the subgraph from Figure \ref{fig:2c_paths_target}. Since $s(v_0 v_1) = +1$, $\varphi(v_1)$ has to be an odd color. Since every odd (resp. even) vertex of the subgraph is only adjacent to even (resp. odd) vertices, we alternate between odd and even colors along the path $\{v_1, v_3, ..., v_{n-1}\}$. Hence $v_{n-1}$ is colored in an even color and completes the homomorphism since $s(v_{n-1} v_0) = -1$.

Suppose not, $G$ is all positive or all negative. If $G$ is an all positive (resp. negative) cycle of odd length, we can color it with the all positive (resp. negative) triangle of $SB$. We have proven that the chromatic number of connected $2$-edge-colored graphs with maximum degree 2 is at most 5.\newline

We will now show that there is no unique graph on 5 vertices that can color the four graphs from Figure \ref{fig:2_counter-example} and therefore that connected $2$-edge-colored graphs with maximum degree 2 are not optimally colorable, that is, there exists a target graph that can color every $2$-edge-colored graph with maximum degree 2.

\begin{figure}[H]
  \begin{center}
    \subfloat[All positive triangle.]{
    	\scalebox{0.8}
		{
			\begin{tikzpicture}[thick]
        		\def \radius {3cm}
        		\def \margin {8} 
        		\tikzstyle{vertex}=[circle,minimum width=1.5em]
        		
        		\node[draw, vertex] (1) at (90:1.2) {};
        		\node[draw, vertex] (2) at (90+120:1.2) {};
        		\node[draw, vertex] (3) at (90-120:1.2) {};
        		
        		\draw (1) -- (2);
        		\draw (2) -- (3);
        		\draw (1) -- (3);
        		\addvmargin{2mm}			
        		\end{tikzpicture}
			
		}
    } \hspace*{0.8cm}
    \subfloat[All negative triangle.]{
    	\scalebox{0.8}
		{
			\begin{tikzpicture}[thick]
        		\def \radius {3cm}
        		\def \margin {8} 
        		\tikzstyle{vertex}=[circle,minimum width=1.5em]
        		
        		\node[draw, vertex] (1) at (90:1.2) {};
        		\node[draw, vertex] (2) at (90+120:1.2) {};
        		\node[draw, vertex] (3) at (90-120:1.2) {};
        		
        		\draw[dashed] (1) -- (2);
        		\draw[dashed] (2) -- (3);
        		\draw[dashed] (1) -- (3);
        		\addvmargin{2mm}
        		\end{tikzpicture}	
			
		}
    } \hspace*{0.8cm}
    \subfloat[Alternating~$C_4$]{
    	\scalebox{0.8}
		{
			\begin{tikzpicture}[thick]
        		\def \radius {3cm}
        		\def \margin {8} 
        		\tikzstyle{vertex}=[circle,minimum width=1.5em]
        		
        		\node[draw, vertex] (1) at (45:1.3) {};
        		\node[draw, vertex] (2) at (45+90:1.3) {};
        		\node[draw, vertex] (3) at (45+180:1.3) {};
        		\node[draw, vertex] (4) at (45-90:1.3) {};
        		
        		\draw[dashed] (1) -- (2);
        		\draw (2) -- (3);
        		\draw[dashed] (3) -- (4);
        		\draw (1) -- (4);
        		\addvmargin{2mm}
        		\end{tikzpicture} 
			
		}
    } \hspace*{0.8cm}
    \subfloat[Alternating~$C_6$]{
    \scalebox{0.7}
		{
			\begin{tikzpicture}[thick]
        		\def \radius {1cm}
        		\def \margin {8} 
        		\tikzstyle{vertex}=[circle,minimum width=1.5em]
        		
        		\node[draw, vertex] (1) at (-60*0+90: 1.2) {$v_1$};
        		\node[draw, vertex] (2) at (-60*1+90: 1.2) {$v_2$};
        		\node[draw, vertex] (3) at (-60*2+90: 1.2) {$v_3$};
        		\node[draw, vertex] (4) at (-60*3+90: 1.2) {$v_4$};
        		\node[draw, vertex] (5) at (-60*4+90: 1.2) {$v_5$};
        		\node[draw, vertex] (6) at (-60*5+90: 1.2) {$v_6$};
        		
        		\draw[dashed] (1) -- (2);
        		\draw (2) -- (3);
        		\draw[dashed] (3) -- (4);
        		\draw (4) -- (5);
        		\draw[dashed] (5) -- (6);
        		\draw (6) -- (1);
        		\addvmargin{2mm}
        		\end{tikzpicture}
			
		}
    }
    \caption{} 
    \label{fig:2_counter-example}  
  \end{center}
\end{figure}

The first three graphs of Figure~\ref{fig:2_counter-example} are 2-edge-colored cliques so they need to be subgraphs of the target graph. There is only one way, up to isomorphisms, to have the two triangles as subgraphs of a 5 vertices graph with a minimal number of edges: the graph $SB$ (Figure~\ref{sub:SB}).

There is only one way, up to isomorphisms, to add edges to $SB$ so that it admits the alternating $C_4$ as a subgraph (see Figure \ref{fig:2_candidate_target}).

\begin{figure}[H]
	\centering
	\scalebox{0.7}
	{
		\begin{tikzpicture}[thick]
		\def \radius {3cm}
		\def \margin {8} 
		\tikzstyle{vertex}=[circle,minimum width=2.5em]
		
		\node[draw, vertex] (5) at (-60*0+90: 2) {$1$};
		\node[draw, vertex] (3) at (-60*1+90: 2) {$2$};
		\node[draw, vertex] (2) at (-60*2+90: 2) {$3$};
		\node[draw, vertex] (4) at (-60*4+90: 2) {$4$};
		\node[draw, vertex] (1) at (-60*5+90: 2) {$5$};
		
		\draw (5) -- (1);
		\draw (5) -- (3);
		\draw (1) -- (3);
		\draw (1) -- (2);
		\draw (3) -- (4);
		
		\draw[dashed] (5) -- (4);
		\draw[dashed] (5) -- (2);
		\draw[dashed] (2) -- (4);
		\draw[dashed] (1) -- (4);
		\draw[dashed] (2) -- (3);
		\end{tikzpicture}
		
	}
	\caption{Candidate target graph on 5 vertices for $2$-edge-colored graphs with maximum degree~2.}
	\label{fig:2_candidate_target}
\end{figure}
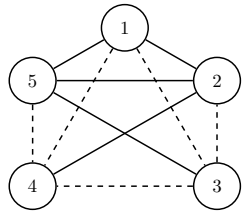

We will now show that it is not possible to color the alternating $C_6$ with our candidate target graph (Figure \ref{fig:2_candidate_target}).
Let $v_1$, $v_2$, ..., $v_6$ be the vertices of the alternating $C_6$ (see Figure \ref{fig:2_counter-example}).
Suppose that color $2$ is used in the coloring. Without loss of generality let $v_1$ be colored in $2$. Since the only negative neighbor of $2$ in the target graph is $3$, $v_2$ needs to be colored in $3$. Similarly, $v_3$ needs to be colored in $5$, $v_4$ in $4$ and $v_5$ in $2$. Since $v_6$ is both a positive and a negative neighbor of vertices colored in $2$, the graph cannot be colored by using the color $2$. 

Similarly, colors $3$, $4$ and $5$ cannot be used to color the alternating $C_6$ so it is not possible to color it with our candidate target graph. We have proven that connected 2-edge-colored graphs are not optimally colorable.
\end{proof}

\subsection{$2$-edge-colored graphs with maximum degree $2$}
\label{sec:2}
While $5$ colors are enough in the case of \emph{connected} $2$-edge-colored graphs with maximum degree $2$ (see Theorem~\ref{thm:2c}), we prove 
in this subsection that $6$ colors are needed when the graphs are not necessarily connected (and this bound is tight).

\begin{theorem}[$\chi_2(\mathcal{D}_2) = 6$]\label{thm:2ecd2}
	The class of $2$-edge-colored graphs with maximum degree $2$ has chromatic number $6$ and is optimally colorable by the target graph depicted in Figure \ref{fig:2_target}.
\end{theorem}

\begin{proof}

The class of graphs with maximum degree 2 is the set of disjoint unions of paths and cycles.

Notice that the graph depicted in Figure \ref{fig:2_target} admits the two graphs from Figure \ref{fig:2c_targets}, $SP_5$ and $SB$, as subgraphs. Therefore, this graph can color any connected 2-edge-colored graph with maximum degree 2 so it can color any 2-edge-colored graph with maximum degree 2.

\begin{figure}[H]
	\centering
	\scalebox{0.8}
	{
		\begin{tikzpicture}[thick]
		\def \radius {3cm}
		\def \margin {8} 
		\tikzstyle{vertex}=[circle,minimum width=1.5em]
		
		\node[draw, vertex] (5) at (-60*0-90: 2) {};
		\node[draw, vertex] (3) at (-60*1-90: 2) {};
		\node[draw, vertex] (2) at (-60*2-90: 2) {};
		\node[draw, vertex] (6) at (-60*3-90: 2) {};
		\node[draw, vertex] (4) at (-60*4-90: 2) {};
		\node[draw, vertex] (1) at (-60*5-90: 2) {};
		
		\draw (5) -- (1);
		\draw (5) -- (3);
		\draw (1) -- (3);
		\draw (1) -- (4);
		\draw (3) -- (2);
		\draw (4) -- (6);
		\draw (2) -- (6);
		
		\draw[dashed] (5) -- (4);
		\draw[dashed] (5) -- (2);
		\draw[dashed] (1) -- (2);
		\draw[dashed] (3) -- (4);
		\draw[dashed] (1) -- (6);
		\draw[dashed] (3) -- (6);
		\draw[dashed] (2) -- (4);
		\end{tikzpicture}
		
	}
	\caption{Target graph for $2$-edge-colored graphs with maximum degree 2.}
	\label{fig:2_target}
\end{figure}
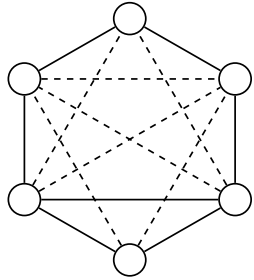

We have shown that the class of graphs with maximum degree 2 has chromatic number at most 6 and is colorable. It is also complete by disjoint union so it is optimally colorable by Lemma~\ref{lem:opt_col}.
By Theorem~\ref{thm:2c}, the class of connected 2-edge-colored graphs with maximum degree 2 has chromatic number 5 and is not optimally colorable. Therefore, there is no single 2-edge-colored graph on 5 vertices that can color every 2-edge-colored paths and cycles. Thus, the class of graphs with maximum degree 2 has chromatic number 6.

\end{proof}
\subsection{Signed graphs with maximum degree $2$}
\label{sec:2s}

In this subsection, we consider the chromatic number of signed graphs with maximum degree~$2$.
Recall that when trying to prove the existence of a homomorphims from a signed graph $G$ to a signed graph $H$, we are allowed to switch a subset of vertices of $G$. This implies that $\chi_s(\mathcal{C}) \le \chi_2(\mathcal{C})$. Therefore, $\chi_s(\mathcal{D}^c_2) \le 5$ by Theorem~\ref{thm:2c} and $\chi_s(\mathcal{D}_2) \le 6$ by Theorem~\ref{thm:2ecd2}. We prove in the following theorem that $4$ colors are enough in both cases (connected or non-connected) and that this is tight.

\begin{theorem}[$\chi_s(\mathcal{D}_2) = 4$]
	The class of signed graphs with maximum degree $2$ has chromatic number $4$ and is optimally colorable by the target graph depicted in Figure \ref{fig:2s_target}.
\end{theorem}

\begin{proof}
An unbalanced $C_4$ is a signed clique by Lemma~\ref{lem:sC} so the chromatic number of signed graphs with maximum degree 2 is at least 4.

\begin{figure}[H]
	\centering
	\scalebox{1}
	{
		\begin{tikzpicture}[thick]
		\def \radius {3cm}
		\def \margin {8} 
		\tikzstyle{vertex}=[circle,minimum width=1.5em]
		
		\node[draw, vertex] (1) at (0, 1.5) {$1$};
		\node[draw, vertex] (2) at (1.5, 1.5) {$2$};
		\node[draw, vertex] (3) at (1.5, 0) {$3$};
		\node[draw, vertex] (4) at (0, 0) {$4$};
		
		\draw (1) -- (2);
		\draw (2) -- (3);
		\draw[dashed] (3) -- (4);
		\draw (1) -- (4);
		\draw (1) -- (3);
		\end{tikzpicture} 
		
	}
	\caption{Target graph $T$ for signed graphs with maximum degree 2.}
	\label{fig:2s_target}
\end{figure}
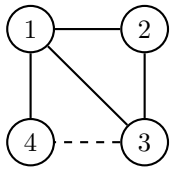

We consider the target graph $T$ depicted in Figure~\ref{fig:2s_target}.

The class of graphs with maximum degree 2 is the set of all paths and cycles. Any signed path is switching equivalent to the all positive path of the same length by Theorem~\ref{thm:Zaslavsky} and every positive path admits a homomorphism to a positive edge. Therefore, a signed path has chromatic number 2.

A cycle of length $n$ is either balanced or unbalanced. 

If it is balanced, it is possible to turn it into an all positive cycle with some number of switches by Theorem \ref{thm:Zaslavsky}. An all positive cycle of even length can be colored with a positive edge. An all positive cycle of odd length can be colored with an all positive triangle (for instance, the subgraph of $T$ induced by vertices $\{1, 2, 3\}$) 

If it is unbalanced, it is possible to turn it into a cycle with exactly one negative edge with some number of switches by Theorem~\ref{thm:Zaslavsky}. Such a cycle of even length can be colored with the cycle $(1, 2, 3, 4)$. Such a cycle of odd length can be colored with the subgraph of $T$ induced by vertices $\{1, 3, 4\}$. See Table \ref{tab:2s_coloring_cases} for reference.

\begin{table}[H]
	\centering
	\begin{tabular}{|>{\centering\arraybackslash}m{.16\linewidth}|>{\centering\arraybackslash}m{.42\linewidth}|>{\centering\arraybackslash}m{.42\linewidth}|N}
		\hline
		& balanced	& unbalanced & \\ \hline		
		even length 
		&  
		
		\begin{tikzpicture}[thick]
		\def \radius {3cm}
		\def \margin {8} 
		\tikzstyle{vertex}=[circle,minimum width=1em]
		
		\node[draw, circle] (1) at (+60*2+60: 1) {};
		\node[draw, circle] (2) at (+60*1+60: 1) {};
		\node[draw, circle] (3) at (    0+60: 1) {};
		\node[draw, circle] (4) at (-60*1+60: 1) {};
		
		\draw (1) -- (2);
		\draw (2) -- (3);
		\draw (3) -- (4);
		\draw[dotted] (1) edge[bend right] (4);		
		
		\node[draw, vertex] (1) at (0+2.25, 1.25-0.125) {$1$};
		\node[draw, vertex] (2) at (1.25+2.25, 1.25-0.125) {$2$};
		
		\draw (1) -- (2);
		
		\draw[->,>=latex] (1.05, 0.5) to (2.05, 0.5);
		\addvmargin{1mm}
		\end{tikzpicture}
		& 
		
		\begin{tikzpicture}[thick]
		\def \radius {3cm}
		\def \margin {8} 
		\tikzstyle{vertex}=[circle,minimum width=1em]
		
		\node[draw, circle] (1) at (+60*2+60: 1) {};
		\node[draw, circle] (2) at (+60*1+60: 1) {};
		\node[draw, circle] (3) at (    0+60: 1) {};
		\node[draw, circle] (4) at (-60*1+60: 1) {};
		
		\draw (1) -- (2);
		\draw[dashed] (2) -- (3);
		\draw (3) -- (4);
		\draw[dotted] (1) edge[bend right] (4);		
		
		\node[draw, vertex] (1) at (0+2.25, 1.25-0.125) {$1$};
		\node[draw, vertex] (2) at (1.25+2.25, 1.25-0.125) {$2$};
		\node[draw, vertex] (3) at (1.25+2.25, 0-0.125) {$3$};
		\node[draw, vertex] (4) at (0+2.25, 0-0.125) {$4$};
		
		\draw (1) -- (2);
		\draw (2) -- (3);
		\draw[dashed] (3) -- (4);
		\draw (1) -- (4);
		
		\draw[->,>=latex] (1.05, 0.5) to (2.05, 0.5);
		\addvmargin{1mm}
		\end{tikzpicture}
		&\\ \hline
		odd length  		
		&
		\begin{tikzpicture}[thick]
		\def \radius {3cm}
		\def \margin {8} 
		\tikzstyle{vertex}=[circle,minimum width=1em]
		
		\node[draw, circle] (1) at (+60*2+60: 1) {};
		\node[draw, circle] (2) at (+60*1+60: 1) {};
		\node[draw, circle] (3) at (    0+60: 1) {};
		\node[draw, circle] (4) at (-60*1+60: 1) {};
		
		\draw (1) -- (2);
		\draw (2) -- (3);
		\draw (3) -- (4);
		\draw[dotted] (1) edge[bend right] (4);		
		
		\node[draw, vertex] (1) at (0+2.25, 1.25-0.125) {$1$};
		\node[draw, vertex] (2) at (1.25+2.25, 1.25-0.125) {$2$};
		\node[draw, vertex] (3) at (1.25+2.25, 0-0.125) {$3$};
		
		\draw (1) -- (2);
		\draw (2) -- (3);
		\draw (1) -- (3);
		
		\draw[->,>=latex] (1.05, 0.5) to (2.05, 0.5);
		\addvmargin{1mm}
		\end{tikzpicture} 
		& 
		\begin{tikzpicture}[thick]
		\def \radius {3cm}
		\def \margin {8} 
		\tikzstyle{vertex}=[circle,minimum width=1em]
		
		\node[draw, circle] (1) at (+60*2+60: 1) {};
		\node[draw, circle] (2) at (+60*1+60: 1) {};
		\node[draw, circle] (3) at (    0+60: 1) {};
		\node[draw, circle] (4) at (-60*1+60: 1) {};
		
		\draw (1) -- (2);
		\draw[dashed] (2) -- (3);
		\draw (3) -- (4);
		\draw[dotted] (1) edge[bend right] (4);		
		
		\node[draw, vertex] (1) at (0+2.25, 1.25-0.125) {$1$};
		\node[draw, vertex] (3) at (1.25+2.25, 0-0.125) {$3$};
		\node[draw, vertex] (4) at (0+2.25, 0-0.125) {$4$};
		
		\draw[dashed] (3) -- (4);
		\draw (1) -- (4);
		\draw (1) -- (3);
		
		\draw[->,>=latex] (1.05, 0.5) to (2.05, 0.5);
		\addvmargin{1mm}
		\end{tikzpicture} &\\ \hline 
	\end{tabular}
	\caption{The subgraph used to color a signed cycle of a given balance and length parity.}
	\label{tab:2s_coloring_cases}
\end{table}
\end{proof}

Note that all balanced cycles of even length have chromatic number 2 (since the target graph needs to contain at least one edge), all cycles of odd length have chromatic number 3 (we need at least 3 vertices since a simple cycle of odd length has chromatic number 3), and all unbalanced cycles of even length have chromatic number 4 (such a cycle contains at least one positive and at least one negative edge so it must also be the case for the target graph and the target graph cannot be a graph on 3 vertices since a cycle of even length cannot admit a homomorphism to a cycle of odd length).

\section{Graphs with maximum degree $k$}\label{sec:dk}
In this section, we present two general theorems that work for any maximum degree $k$. The first one requires us to first find for each $k$ a target graph that has some special properties while the second one gives us directly an upper bound for every $k$ (at the cost of giving a looser upper bound). \newline

A graph is said to be \textit{$k$-degenerate} if each of its subgraphs contains at least one vertex of degree at most $k$.
\begin{lemma}
\label{lem:k-deg}
If $T$ is a 2-edge-colored graph with Property $P_{k-1, \floor*{\frac{k-1}{2}}+1}$, then every $(k-1)$-degenerate 2-edge-colored graph with maximum degre $k$ admits a homomorphism to $T$.
\end{lemma}

\begin{proof}
Let $T$ be a 2-edge-colored graph with Property $P_{k-1, \floor*{\frac{k-1}{2}}+1}$. Let $G$ be a $(k-1)$-degenerate 2-edge-colored graph with maximum degree $k$. We proceed by induction on the number of vertices of $G$. Let $s$ be the signature of $G$. Let $u \in V(G)$ be a vertex of degree $k-1$, $v_1, v_2, ..., v_{l-1}, v_l, v_{l+1}, ..., v_{k-1}$ be its neighbors such that $s(v_1 u) = s(v_2 u) = ... = s(v_l u) \neq s(v_{l+1} u) = s(v_{l+2} u) = ... = s(v_{k-1})$ and $l \leq \floor*{\frac{k-1}{2}}$. 

By the induction hypothesis, $G-u$ admits a homomorphism $\varphi$ to $T$.
By Property $P_{k-1, \floor*{\frac{k-1}{2}}+1}$ of $T$, for each $v_{l+i}$ ($1 \leq i \leq k-1-l$) we can recolor $v_{l+i}$ (if needed) such that $\varphi(v_{l+i}) \neq \varphi(v_{j})$ for $1 \leq j \leq l$.

We can now use $P_{k-1, \floor*{\frac{k-1}{2}}+1}$ of $T$ to extend $\varphi$ to $G$ (i.e. to color $u$). Remember that if two of the $v_i$ have the same color and are adjacent to $u$ with same sign they induce the same constraints on the coloring of $u$ and do no prevent us from using $P_{k-1, \floor*{\frac{k-1}{2}}+1}$.
\end{proof}

\begin{lemma}
\label{lem:+2}
If all the $(k-1)$-degenerate 2-edge-colored graphs with maximum degree k admit a homomorphism to a single edge-transitive target graph on $n$ vertices then all the graphs in $\mathcal{D}_k$ admit a homomorphism to a single target graph on $n+2$ vertices.
\end{lemma}

\begin{proof}
Let $T$ be an edge-transitive target graph on $n$ vertices that can color every $(k-1)$-degenerate 2-edge-colored graph with maximum degree $k$. Let $xy$ be a positive edge of $T$. Consider the graph $T^*$ obtained from $T$ by adding two new vertices $x'$ and $y'$ as follows. Link $x'$ and $y'$ to the vertices of $T$ in the same way as $x$ and $y$ are, respectively; add an edge $x'y'$ with $s(x'y') = -1$; finally we add edges $xx'$ and $yy'$ with $s(xx') = -1$ and $s(yy') = +1$. To prove that every graph from $\mathcal{D}_k$ admits a homomorphism to $T^*$ it suffices to show that every connected $k$-regular graph admits a homomorphism to $T^*$. \newline

Let $G$ be a $k$-regular 2-edge-colored graph. Since $T$ can color every $(k-1)$-degenerate graph with maximum degree $k$, $T$ contains an all positive $K_{k}$ as a subgraph. Since $T$ is edge-transitive, it is in particular vertex-transitive and there exists an all positive $K_{k}$ $\{y, y_1, y_2, ..., y_{k-1}\}$ that contains $y$. Since $y$ and $y'$ have the same neighborhoods in $T$ and they are adjacent with a positive edge, $\{y, y', y_1, y_2, ..., y_{k-1}\}$ is an all positive $K_{k+1}$. If $G$ is all positive, it can be colored using this all positive $K_{k+1}$. We can now assume that $G$ contains at least one negative edge.

Let $uv$ be a negative edge. The graph $G-uv$ is $(k-1)$-degenerate graph and admits a homomorphism $\varphi$ to $T$.

If $\varphi(u) \varphi(v)$ is a negative edge in $T$ then $\varphi$ is already a homomorphism from $G$ to $T^*$. 

If $\varphi(u) \varphi(v)$ is a positive edge in $T$, then by the edge-transitivity of $T$ there exists a homomorphism $\varphi'$ from $G$ to $T$ such that $\varphi'(u) = x$ and $\varphi'(v) = y$. The following application $\varphi''$ is a homomorphism from $G$ to $T^*$ because $x' y'$ is a negative edge and $x'$ and $y'$ have the same positive and negative neighbors in $T^*$ as $x$ and $y$ in $T$.

$$
\varphi''(w) = \left\{
    \begin{array}{ll}
        x' & \mbox{if } w = u, \\
        y' & \mbox{if } w = v, \\
        \varphi'(w) & \mbox{otherwise.}\\
    \end{array}
\right.
$$

If $\varphi(u) = \varphi(v)$, then by the vertex-transitivity of $T$ there exists a homomorphism $\varphi'$ from $G$ to $T$ such that $\varphi'(u) = \varphi'(v) = x$. The following application $\varphi''$ is a homomorphism from $G$ to $T^{*}$ because $x' x$ is a negative edge and $x'$ has the same positive and negative neighbors in $T^*$ as $x$ in $T$.

$$
\varphi''(w) = \left\{
    \begin{array}{ll}
        x' & \mbox{if } w = u, \\
        \varphi'(w) & \mbox{otherwise.}\\
    \end{array}
\right.
$$

\end{proof}

\begin{theorem}
\label{thm:META}
If there exists an edge-transitive 2-edge-colored graph $T$ with Property $P_{k-1, \floor*{\frac{k-1}{2}}+1}$ on $n$ vertices then the class of (connected) 2-edge-colored graphs with maximum degree $k$ has chromatic number at most $n+2$.
\end{theorem}

\begin{proof}
This follows from Lemmas~\ref{lem:k-deg} and \ref{lem:+2}.
\end{proof}
    
\begin{corollary}
\label{cor:META}
If there exists an edge-transitive antitwinned 2-edge-colored graph $\rho(T)$ with Property $P_{k-1, \floor*{\frac{k-1}{2}}+1}$ on $2n$ vertices then the class of (connected) signed graphs with maximum degree $k$ has chromatic number at most $n + 2$.
\end{corollary}

\begin{proof}
By Lemma~\ref{lem:k-deg}, $\rho(T)$ can color all $(k-1)$-degenerate 2-edge-colored graphs with maximum degree $k$ (and by Lemma~\ref{lem:BG}, $T$ can color all $(k-1)$-degenerate signed graphs with maximum degree $k$).  

We apply Lemma~\ref{lem:+2} to get a target graph $\rho(T)^*$ that can color every graph in $\mathcal{D}_k$. This graph is not antitwinned since $u$ and $v$ do not have antitwins. We add the missing antitwins of $u'$ and $v'$ in order to get an antitwinned signed target graph on $2n + 4$ vertices. By Lemma~\ref{lem:BG} we get that every signed graph in $\mathcal{D}_k$ admits a homomorphism to a single target graph on $n+2$ vertices.
\end{proof}

By Lemma~\ref{lem:PTRSP}, $TR(SP_5)$ is edge-transitive and has Property $P_{2,2}$, $TR(SP_{13})$ is edge-transitive and has Property $P_{3, 2}$, and $TR(SP_{53})$ is edge-transitive and has Property $P_{4, 4}$ (calculated by computer). Therefore, we can deduce the following bounds using Theorem~\ref{thm:META} and Corollary~\ref{cor:META}:
\begin{itemize}
\item $8\le \chi_2(\mathcal{D}^c_3) \le \chi_2(\mathcal{D}_3) \le 14$ (by Theorems~\ref{thm:2ecc} and~\ref{thm:META} using $TR(SP_5)$);
\item $\chi_s(\mathcal{D}^c_3) \le \chi_s(\mathcal{D}_3) \le 8$ (by Corollary~\ref{cor:META} using $TR(SP_5)$);
\item $12\le \chi_2(\mathcal{D}^c_4) \le \chi_2(\mathcal{D}_4) \le 30$ (by Theorems~\ref{thm:2ecc} and~\ref{thm:META} using $TR(SP_{13})$);
\item $10\le \chi_s(\mathcal{D}^c_4) \le \chi_s(\mathcal{D}_4) \le 16$ (by Theorem~\ref{thm:sc} and by Corollary~\ref{cor:META} using $TR(SP_{13})$);
\item $16\le \chi_2(\mathcal{D}^c_5) \le \chi_2(\mathcal{D}_5) \le 110$ (by Theorems~\ref{thm:2ecc} and~\ref{thm:META} using $TR(SP_{53})$);
\item $12\le \chi_s(\mathcal{D}^c_5) \le \chi_s(\mathcal{D}_5) \le 56$ (by Theorem~\ref{thm:sc} and Corollary~\ref{cor:META} using $TR(SP_{53})$);
\end{itemize}

We present in Section~\ref{sec:3} two theorems that yield better results for 2-edge-colored graphs with maximum degree~$3$, namely the above-mentioned upper bound of $14$ is decreased to $11$ for non-necessarily connected graphs (see Subsection~\ref{sec:3_SP9+2}) and to $10$ for connected ones (see Subsection~\ref{sec:3c}). The chromatic number of signed graphs with maximum degree~3 can also be improved (see Subsection~\ref{sec:3s}).\newline

We could find other upper bounds for the chromatic number of 2-edge-colored and signed graphs with maximum degree $k\geq 5$ using Theorem~\ref{thm:META} and Corollary~\ref{cor:META} by calculating the properties of $TR(SP_q)$ and $\rho(SP_q)$ for greater values of $q$. This would require a lot of processing time or a better way to compute or approximate these properties. \newline

We now present a general upper bound for the chromatic number of 2-edge-colored and signed graphs with maximum degree $k$ that does not require computations.

An \textit{$(m, n)$-colored-mixed graph} is a graph in which each pair of vertices can either be connected by an edge, of which there are $n$ types (in the same way there are 2 types of edges in a 2-edge-colored graph) or an arc (an edge with an orientation represented by an ordered pair of vertices instead of a 2-set), of which there are $m$ types. 

A 2-edge-colored graph is therefore a $(0, 2)$-colored-mixed graph.

Das, Nandi and Sen proved the following general theorem on $(m, n)$-colored-mixed graphs using a probabilistic argument.

\begin{theorem}[\cite{MIXED}]
\label{thm:mixed}
The chromatic number of a connected $(m, n)$-colored-mixed graph with maximum degree $k \geq 5$ is at most $2 \cdot (k-1)^{2m+n} \cdot (2m+n)^{k-1}+2$ and at least $(2m+n)^{\frac{k}{2}}$.
\end{theorem}

From this general theorem on colored-mixed graphs we can obtain the following corollary on 2-edge-colored graphs.

\begin{corollary}[$2^{\frac{k}{2}} \leq \chi_2(\mathcal{D}_k^c) \leq (k-1)^2 \cdot 2^{k} + 2$]
\label{cor:k}
The chromatic number of connected 2-edge-colored graphs with maximum degree $k \geq 5$ is at most $(k-1)^2 \cdot 2^{k} + 2$ and at least $2^{\frac{k}{2}}$.
\end{corollary}

The upper bound also applies trivially to connected signed graphs. The lower bound given by Theorem~\ref{thm:2ecc} is better than $2^{\frac{k}{2}}$ for $k \leq 10$.

\begin{corollary}[$2^{\frac{k}{2}} \leq \chi_2(\mathcal{D}_k) \leq k^2 \cdot 2^{k+1}$]\label{cor:upper}
The chromatic number of 2-edge-colored graphs with maximum degree $k \geq 5$ is at most $k^2 \cdot 2^{k+1}$ and at least $2^{\frac{k}{2}}$.
\end{corollary}
\begin{proof}
The lower bound from Corollary~\ref{cor:k} also applies trivially to disconnected graphs.

Theorem~\ref{thm:mixed} is proved by showing that there exist an $(m, n)$-colored-mixed graph with Property $Q^{k-1, j}_{1+(k-j)(k-2)}$ on $2(k-1)^p\cdot p^{k-1}+2$ vertices where $p = 2m+n \geq 2$ and $k\geq 5$. We will not give the definition of properties of the type $Q^{t, j}_{g(j)}$ here but for the case of 2-edge-colored graphs, Property $Q^{k-1, j}_{1+(k-j)(k-2)}$ implies Property $P_{k-1, k-1}$.

Therefore, there exists a 2-edge-colored graph with Property $P_{k, k}$ on $k^2\cdot 2^{k+1}$ vertices. We conclude using Lemma~\ref{lem:k-deg} (a $k$-regular graph is also $k$-degenerate).
\end{proof}

The upper bound given in Corollary~\ref{cor:upper} also applies trivially to signed graphs. The following theorem gives a lower bound for the chromatic number of signed graphs with maximum degree $k \geq 5$.

\begin{theorem}[$2^{\frac{k}{2} - 1} \leq \chi_s(\mathcal{D}_k)$]
The chromatic number of signed graphs with maximum degree $k \geq 5$ is at least $2^{\frac{k}{2} - 1}$.
\end{theorem}

\begin{proof}
We adapt the proof of the lower bound of Theorem~\ref{thm:mixed} for signed graphs.

Let $G$ be a labeled connected simple graph. We denote by $\chi_s(G)$ the maximum of the chromatic numbers of all the signed graphs with underlying graph $G$.

The number of labeled signed graphs with underlying graph $G$ is $2^{|E(G)|}$ since each edge of $G$ can either be positive or negative.

For each of these signed graphs, there are $2^{|V(G)|-1}$ ways to switch its vertices (note that switching all and none of the vertices yields the same signed graph).

Each of these signed graphs has chromatic number at most $\chi_s(G)$ so it admits a homomorphism to at least one complete signed graph on $\chi_s(G)$ vertices. There are $2^{{\chi_s(G)}\choose{2}}$ complete labeled signed graphs on $\chi_s(G)$ vertices.

There are $\chi_s(G)^{|V(G)|}$ applications from the vertex set of a graph on $|V(G)|$ vertices to the vertex set of a graph on $\chi_s(G)$ vertices.

For each of the labeled signed graphs with underlying graph $G$, for at least one of its switching equivalent graphs, at least one of the applications from the vertex set of this graph to the vertex set of at least one of the complete signed graphs on $\chi_s(G)$ vertices is a homomorphism. Therefore we have:

$$ 2^{|V(G)|-1} \cdot \chi_s(G)^{|V(G)|} \cdot 2^{{\chi_s(G)}\choose{2}} \geq 2^{|E(G)|}$$

Remark: Let $G^1$, $G^2$ ($G^1 \neq G^2$) be two of the $2^{|E(G)|}$ labeled signed graphs with underlying graph $G$. Graphs $G^1$ and $G^2$ have a signature that is different on a least one edge and therefore an application from the vertex set of $G$ to a given complete signed graph on $\chi_s(G)$ after switching the same subset of vertices in $G^1$ and $G^2$ cannot be a homomorphism for both $G^1$ and $G^2$.

We raise each side to $\frac{1}{|V(G)|}$:

$$\displaystyle 2^{\frac{|V(G)|-1}{|V(G)|}} \cdot \chi_s(G)^{\frac{|V(G)|}{|V(G)|}} \cdot 2^{{\frac{{{\chi_s(G)}\choose{2}}}{|V(G)|}}} \geq 2^{\frac{|E(G)|}{|V(G)|}}$$

$$ \chi_s(G) \geq \frac{ 2^{\frac{|E(G)|}{|V(G)|}} }{ 2^{\frac{|V(G)|-1}{|V(G)|}} \cdot 2^{{\frac{{{\chi_s(G)}\choose{2}}}{|V(G)|}}}}$$

We choose $G$ $k$-regular:

$$ \chi_s(G) \geq \frac{ 2^{\frac{k}{2}} }{ 2^{\frac{|V(G)|-1}{|V(G)|}} \cdot 2^{{\frac{{{\chi_s(G)}\choose{2}}}{|V(G)|}}}}$$

Since $\chi_s(G)$ is bounded (by Corollary~\ref{cor:k} and the fact that $\chi_s(G) \leq \chi_2(G)$ for any $G$), the right side approaches $2^{\frac{k}{2}-1}$ as $|V(G)|$ goes to infinity.

\end{proof}

\section{Graphs with maximum degree 3}
In this section, we consider graphs with maximum degree~3 and we improve the upper bounds that were obtained by Theorem~\ref{thm:META} and Corollary~\ref{cor:META} in Section~\ref{sec:dk}.
\label{sec:3}
\subsection{$2$-edge-colored graphs with maximum degree $3$}
\label{sec:3_SP9+2}

\medmuskip=0mu

By Theorem~\ref{thm:2ecc}, the class of 2-edge-colored graphs with maximum degree 3 has chromatic number at least 8.

The graph $TR(SP_5)$ is the smallest known edge-transitive 2-edge-colored graph with Property $P_{2, 2}$. By Theorem~\ref{thm:META}, every 2-edge-colored graph with maximum degree 3 has chromatic number at most $|TR(SP_5)|+2 = 14$.
In this subsection, we adapt the proof to work with the target graph $SP_9$ even though it does not have property $P_{2, 2}$, decreasing the upper bound to $|SP_9|+2=11$.\newline

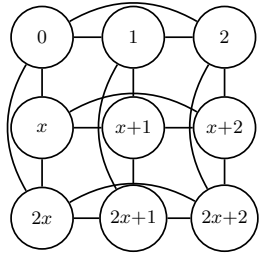
\begin{figure}[H]
	\centering
	\scalebox{0.8}
	{
		\begin{tikzpicture}[thick]
		\def \radius {1cm}
		\def \margin {8} 
		\tikzstyle{vertex}=[circle,minimum width=2.9em]
		
		\node[draw, vertex] (0) at (0, 3) {\small$0$};
		\node[draw, vertex] (1) at (1.5, 3) {\small$1$};
		\node[draw, vertex] (2) at (3, 3) {\small$2$};
		\node[draw, vertex] (3) at (0, 1.5) {\small$x$};
		\node[draw, vertex] (4) at (1.5, 1.5) {\small$x+1$};
		\node[draw, vertex] (5) at (3, 1.5) {\small$x+2$};
		\node[draw, vertex] (6) at (0, 0) {\small$2x$};
		\node[draw, vertex] (7) at (1.5, 0) {\small$2x+1$};
		\node[draw, vertex] (8) at (3, 0) {\small$2x+2$};

		\draw (0) -- (1);
		\draw (1) -- (2);
		\draw (2) edge[bend right]  (0);
		
		\draw (3) -- (4);
		\draw (4) -- (5);
		\draw (5) edge[bend right]  (3);
		
		\draw (6) -- (7);
		\draw (7) -- (8);
		\draw (8) edge[bend right]  (6);
		
		\draw (0) -- (3);
		\draw (3) -- (6);
		\draw (6) edge[bend left]  (0);
		
		\draw (1) -- (4);
		\draw (4) -- (7);
		\draw (7) edge[bend left]  (1);
		
		\draw (2) -- (5);
		\draw (5) -- (8);
		\draw (8) edge[bend left]  (2);
		
		\end{tikzpicture}
	}
	\caption{The graph $SP_9$, non-edges are negative edges.}
	\label{fig:SP_9}
\end{figure}

Graph $SP_9$ (see Figure~\ref{fig:SP_9}) has Properties $P_{1, 4}$ and $P_{2,1}$ by Lemma~\ref{lem:PSP}. We also introduce the following new property of $SP_9$.

\begin{lemma}[Property $P_{2, 2}^*$ of $SP_9$]
\label{P22star}
Given two vertices $u$ and $v$ of $SP_9$ and two signs $\{s_1, s_2\} \in \{-1, +1\}^2$ such that $|\{s(uv), s_1, s_2\}| > 1$, there are two vertices $w_1$ and $w_2$ of $SP_9$ such that $s(uw_1)=s(uw_2)=s_1$ and $s(vw_1)=s(vw_2)=s_2$.
\end{lemma}
\begin{proof}
Since $SP_9$ is edge-transitive and antiautomorphic by Lemma~\ref{lem:PSP}, it suffices to consider the case $u=0$ and $v=1$. Since $01$ is a positive edge, we have two cases to consider:
\begin{itemize}
\item Either $s_1 = s_2 = -1$ and we can have $w_1 = x+2$ and $w_2 = 2x+2$;
\item Or $s_1 = +1$, $s_2 = -1$ and we can have $w_1 = x$ and $w_2 = 2x$.
\end{itemize}
\end{proof}

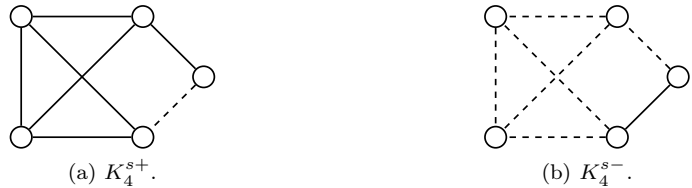
\begin{figure}[H]
  \begin{center}
    \subfloat[$K_4^{s+}$.]{
    	\scalebox{0.8}
		{
			\begin{tikzpicture}[thick][thick]
        		\def \radius {1cm}
        		\def \margin {8} 
        		\tikzstyle{vertex}=[circle,minimum width=1em]
        		
        		\node[draw, vertex] (1) at (0,0) {};
        		\node[draw, vertex] (2) at (-1,-1) {};
        		\node[draw, vertex] (3) at (-1,1) {};
        		\node[draw, vertex] (6) at (-3,1) {};
        		\node[draw, vertex] (7) at (-3,-1) {};
        		
        		\draw[dashed] (1) -- (2);
        		\draw (1) -- (3);
        		\draw  (2) -- (6);
        		\draw  (2) -- (7);
        		\draw (3) -- (6);
        		\draw (3) -- (7);
        		\draw (6) -- (7);
        		\end{tikzpicture}
			
		}
    } \hspace*{3cm}
    \subfloat[$K_4^{s-}$.]{
    \scalebox{0.8}
		{
			\begin{tikzpicture}[thick][thick]
        		\def \radius {1cm}
        		\def \margin {8} 
        		\tikzstyle{vertex}=[circle,minimum width=1em]
        		
        		\node[draw, vertex] (1) at (0,0) {};
        		\node[draw, vertex] (2) at (-1,-1) {};
        		\node[draw, vertex] (3) at (-1,1) {};
        		\node[draw, vertex] (6) at (-3,1) {};
        		\node[draw, vertex] (7) at (-3,-1) {};
        		
        		\draw (1) -- (2);
        		\draw[dashed] (1) -- (3);
        		\draw[dashed] (2) -- (6);
        		\draw[dashed]  (2) -- (7);
        		\draw[dashed] (3) -- (6);
        		\draw[dashed] (3) -- (7);
        		\draw[dashed] (6) -- (7);
        		\end{tikzpicture}
			
		}
    }
    \caption{The two $K_4^s$ graphs.} 
    \label{fig:3_K4s}  
  \end{center}
\end{figure}

We say that a 2-edge-colored graph is a $K_4^{s}$ if it is the graph $K_4^{s+}$ or $K_4^{s-}$, the all positive or the all negative complete graph on $4$ vertices with one edge subdivided into a path of length 2 with one negative and one positive edge. See Figure~\ref{fig:3_K4s}.

\begin{lemma}
\label{lem:3_2deg}
Every 2-degenerate 2-edge-colored graph with maximum degree 3 that does not contain a $K_4^s$ as a subgraph admits a homomorphism to $SP_9$. 
\end{lemma}

\begin{proof}
We prove the lemma by induction on the number of vertices.

Let $G$ be a 2-degenerate 2-edge-colored graph with maximum degree 3 that does not contain a $K_4^s$ as a subgraph. Let $s$ be the signature of $G$. 

Suppose that $G$ contains a vertex $u$ of degree 1. By the induction hypothesis, $G-u$ admits a homomorphism $\varphi$ to $SP_9$. By Property $P_{1, 4}$ of $SP_9$, we can extend this homomorphism to $G$. We thus assume that $G$ does not contain a vertex of degree $1$.

Let $u \in V(G)$ be a vertex of degree 2 and $v$ and $w$ be its neighbors.
By the induction hypothesis, $G-u$ admits a homomorphism $\varphi$ to $SP_9$.

Suppose that $s(vu) = s(wu)$. If $\varphi(v) \neq \varphi(w)$, it is possible to extend $\varphi$ to $G$ (i.e. to color $u$) by Property $P_{2, 1}$ of $SP_9$. If $\varphi(u) = \varphi(v)$, it is still possible to extend $\varphi$ to $G$ because $v$ and $w$ give the same constraints on the color of $u$. We thus assume that $s(vu) \neq s(wu)$

We show in the remainder of this proof that it is always possible to recolor $G-u$ such that $v$ and $w$ get distinct colors (if it is not already the case). Once $\varphi(v) \neq \varphi(w)$, $\varphi$ can be extended to $G$ by Property $P_{2, 1}$ of $SP_9$. \newline

If $G - u$ has two components (if $u$ is a cut-vertex) we can apply the induction hypothesis to both components. By vertex-transitivity of $SP_9$ we can recolor $w$ such that $\varphi(v) \neq \varphi(w)$. We thus assume that $G-u$ is connected.

If $v$ and $w$ are adjacent, $\varphi(v) \neq \varphi(w)$. We thus assume that $v$ and $w$ are not adjacent.

If $v$ has degree 2 in $G$, Property $P_{1, 4}$ of $SP_9$ ensures that we can recolor $v$ with a color distinct from $\varphi(w)$. We thus assume that $v$ and $w$ have degree 3.

Let $v_1$ and $v_2$ be the other two neighbors of $v$ and let $w_1$ and $w_2$ be the other two neighbors of $w$.

Suppose that $w_1$ and $w_2$ are not adjacent. Let $G'$ be $G-\{u, w\}$ with an added edge between $w_1$ and $w_2$ that does not have the same sign as $w_1 w$. The graph $G'$ is 2-degenerate (since we assumed that $G-u$ is connected) so it admits a homomorphism to $SP_9$ by the induction hypothesis. Since $s(w_1 w_2) \neq s(w_1 w)$, we can now apply Property $P_{2, 2}^*$ of $SP_9$ to extend this homomorphism to $G-u$ such that $v$ and $w$ get distinct colors. We thus assume that $w_1 w_2 \in E(G)$.

If $|\{s(w_1 w_2), s(w_1 w), s(w_2 w)\}| > 1$, then Property $P_{2, 2}^*$ of $SP_9$ applies and ensures that we can recolor $w$ with a color distinct from $\varphi(v)$. We thus assume that $s(w_1 w_2) = s(w_1 w) = s(w_2 w)$.

Suppose that $v$ and $w_1$ are not adjacent. If $w_1 w$ is a positive (resp. negative) edge, let $G'$ be $G$ after removing $u$ and $w$ and adding a negative (resp. positive) edge between $v$ and $w_1$. By the induction hypothesis, $G'$ admits a homomorphism $\varphi'$ to $SP_9$. Since $w_1$ and $w_2$ are adjacent, we necessarily have $\varphi'(w_1) \neq \varphi'(w_2)$ and thus Property $P_{2, 1}$ of $SP_9$ allows us to extend $\varphi'$ to $G - u$ (i.e. to give a color to $w$). Since $v$ is a negative (resp. positive) neighbor of $w_1$ and $w$ is a positive (resp. negative) neighbor of $w_1$, $v$ and $w$ have distinct colors. We can now assume that $w_1$ is adjacent to $v$ and by symmetry we can also assume that $w_2$ is adjacent to $v$.

Since we assumed that $s(w_1 w_2) = s(w_1 v) = s(w_1 w) = s(w_2 v) = s(w_2 w)$ and $s(v u) \neq s(w u)$, $G$ contains a $K_4^s$ as a subgraph, a contradiction.
\end{proof}

Consider the graph $SP^*_9$ obtained from $SP_9$ by adding two new vertices $0'$ and $1'$ as follows. Take the two vertices $0$ and $1$ of $SP_9$ (note that $s(01) = +1$), and link $0'$ and $1'$ to the vertices of $SP_9$ in the same way as $0$ and $1$ are, respectively; add an edge $0'1'$ with $s(0'1') = -1$; finally we add edges $00'$ and $11'$ with $s(00') = -1$ and $s(11') = +1$.

\begin{lemma}
\label{lem:3_3reg}
Every 3-regular 2-edge-colored graph that does not contain a $K_4^s$ admits a homomorphism to $SP_9^{*}$.
\end{lemma}
\begin{proof}
We follow the proof of Lemma~\ref{lem:+2} while using Lemma~\ref{lem:3_2deg} instead of Lemma~\ref{lem:k-deg}.
\end{proof}

\begin{theorem}[$\chi_2(\mathcal{D}_3) \leq 11$]
\label{thm:3}
	The class of $2$-edge-colored graphs with maximum degree $3$ has chromatic number at most 11 and is optimally colorable by $SP_9^{*}$.
\end{theorem}
\begin{proof}

Let $G$ be a 2-edge-colored graph with maximum degree 3 and let $C$ be a component of $G$. It suffices to prove that every component admits a homomorphism to $SP_9^{*}$. \newline

Suppose $C$ is 2-degenerate or contains a $K_4^s$ as a subgraph. Let $C'$ be obtained from $C$ after removing all its $K_4^s$. The component $C'$ is thus 2-degenerate with maximum degree 3 and does not contain a $K_4^s$. By Lemma~\ref{lem:3_2deg}, $C'$ admits a homomorphism $\varphi$ to $SP_9$. Figure~\ref{fig:3_stitch} shows how to extend $\varphi$ to $C$ and $SP_9^*$ for a $K_4^{s+}$ if the edge linking it to the rest of the graph is positive. There are 4 cases depending on the color of the vertex $v$ that links a $K_4^{s+}$ to the rest of the graph. Since $0'$ (resp. $1'$) has the same neighboorhood as $0$ (resp $1$) in $SP_9$, $SP_9$ is antiautomorphic, and $K_4^{s+}$ is antiisomorphic to $K_4^{s-}$ (it is isomorphic to $K_4^{s-}$ after replacing each positive edge by a negative one and vice versa), this can also be done for a $K_4^{s-}$ or if the edge linking a $K_4^{s}$ to the rest of the graph is negative. 

\begin{figure}[H]
  \begin{center}
    \subfloat[$\varphi(v) \in \{0, x+2, 2x$\}]{
    	\scalebox{0.7}
		{
			\begin{tikzpicture}[thick][thick]
        		\def \radius {1cm}
        		\def \margin {8} 
        		\tikzstyle{vertex}=[circle,minimum width=1em]
        		
        		\node[draw, vertex] (1) at (0,0) {};
        		\node[draw, vertex] (2) at (-1,-1) {};
        		\node[draw, vertex] (3) at (-1,1) {};
        		\node[draw, vertex] (6) at (-3,1) {};
        		\node[draw, vertex] (7) at (-3,-1) {};
        		\node[draw, vertex] (8) at (1, 0) {$v$};
        		
        		\draw (1.north) node[above]{$x$};
        		\draw (2.north) node[below, shift={(0,-0.3)}]{$1$};
        		\draw (3.north) node[above]{$x+1$};
        		\draw (6.north) node[above]{$1'$};
        		\draw (7.north) node[below, shift={(0,-0.3)}]{$2x+1$};
        		
        		\draw[dashed] (1) -- (2);
        		\draw (1) -- (3);
        		\draw  (2) -- (6);
        		\draw  (2) -- (7);
        		\draw (3) -- (6);
        		\draw (3) -- (7);
        		\draw (6) -- (7);
        		\draw (1) -- (8);
        		\end{tikzpicture}
			
		}
    } 
    \subfloat[$\varphi(v) \in \{0, x+2, 2x$\}]{
    	\scalebox{0.7}
		{
			\begin{tikzpicture}[thick][thick]
        		\def \radius {1cm}
        		\def \margin {8} 
        		\tikzstyle{vertex}=[circle,minimum width=1em]
        		
        		\node[draw, vertex] (1) at (0,0) {};
        		\node[draw, vertex] (2) at (-1,-1) {};
        		\node[draw, vertex] (3) at (-1,1) {};
        		\node[draw, vertex] (6) at (-3,1) {};
        		\node[draw, vertex] (7) at (-3,-1) {};
        		\node[draw, vertex] (8) at (1, 0) {$v$};
        		
        		\draw (1.north) node[above]{$x+2$};
        		\draw (2.north) node[below, shift={(0,-0.3)}]{$1$};
        		\draw (3.north) node[above]{$x+1$};
        		\draw (6.north) node[above]{$1'$};
        		\draw (7.north) node[below, shift={(0,-0.3)}]{$2x+1$};
        		
        		\draw[dashed] (1) -- (2);
        		\draw (1) -- (3);
        		\draw  (2) -- (6);
        		\draw  (2) -- (7);
        		\draw (3) -- (6);
        		\draw (3) -- (7);
        		\draw (6) -- (7);
        		\draw (1) -- (8);
        		\end{tikzpicture}
			
		}
    }
    \subfloat[$\varphi(v) = 2x+1$]{
    	\scalebox{0.7}
		{
			\begin{tikzpicture}[thick][thick]
        		\def \radius {1cm}
        		\def \margin {8} 
        		\tikzstyle{vertex}=[circle,minimum width=1em]
        		
        		\node[draw, vertex] (1) at (0,0) {};
        		\node[draw, vertex] (2) at (-1,-1) {};
        		\node[draw, vertex] (3) at (-1,1) {};
        		\node[draw, vertex] (6) at (-3,1) {};
        		\node[draw, vertex] (7) at (-3,-1) {};
        		\node[draw, vertex] (8) at (1, 0) {$v$};
        		
        		\draw (1.north) node[above]{$2x$};
        		\draw (2.north) node[below, shift={(0,-0.3)}]{$1$};
        		\draw (3.north) node[above]{$2x+1$};
        		\draw (6.north) node[above]{$1'$};
        		\draw (7.north) node[below, shift={(0,-0.3)}]{$x+1$};
        		
        		\draw[dashed] (1) -- (2);
        		\draw (1) -- (3);
        		\draw  (2) -- (6);
        		\draw  (2) -- (7);
        		\draw (3) -- (6);
        		\draw (3) -- (7);
        		\draw (6) -- (7);
        		\draw (1) -- (8);
        		\end{tikzpicture}
			
		}
    } 
    \subfloat[$\varphi(v) = 1$]{
    \scalebox{0.7}
		{
			\begin{tikzpicture}[thick][thick]
        		\def \radius {1cm}
        		\def \margin {8} 
        		\tikzstyle{vertex}=[circle,minimum width=1em]
        		
        		\node[draw, vertex] (1) at (0,0) {};
        		\node[draw, vertex] (2) at (-1,-1) {};
        		\node[draw, vertex] (3) at (-1,1) {};
        		\node[draw, vertex] (6) at (-3,1) {};
        		\node[draw, vertex] (7) at (-3,-1) {};
        		\node[draw, vertex] (8) at (1, 0) {$v$};
        		
        		\draw (1.north) node[above]{$0$};
        		\draw (2.north) node[below, shift={(0,-0.3)}]{$2$};
        		\draw (3.north) node[above]{$1$};
        		\draw (6.north) node[above]{$0$};
        		\draw (7.north) node[below, shift={(0,-0.3)}]{$1'$};
        		
        		\draw[dashed] (1) -- (2);
        		\draw (1) -- (3);
        		\draw  (2) -- (6);
        		\draw  (2) -- (7);
        		\draw (3) -- (6);
        		\draw (3) -- (7);
        		\draw (6) -- (7);
        		\draw (1) -- (8);
        		\end{tikzpicture}
			
		}
    }
    \caption{How to extend $\varphi$ to a $K_4^{s+}$ depending on the color of the vertex $v$ that links the $K_4^{s+}$ to the rest of the graph.} 
    \label{fig:3_stitch}
  \end{center}
\end{figure}
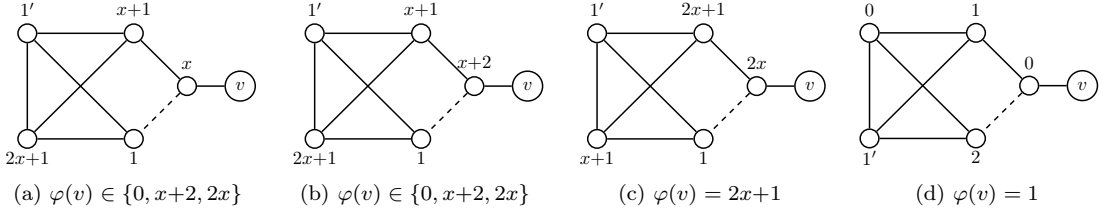

If $C$ is not 2-degenerate and does not contain $K_4^s$ as a subgraph, we can conclude by using Lemma~\ref{lem:3_3reg}.
\end{proof}

\medmuskip=4mu
\subsection{Connected 2-edge-colored graphs with maximum degree 3}\label{sec:3c}
In this subection, we consider the \emph{connected} 2-edge-colored graphs with maximum degree 3. In the previous subsection, we proved that $11$ colors are enough (when the graphs are not necessarily connected) by proving the existence of a universal target graph on $11$ vertices $SP_9^*$. In the connected case, we decrease the upper bound to $10$ by using multiple target graphs on $10$ vertices.

By Theorem~\ref{thm:2ecc}, the class of connected 2-edge-colored graphs with maximum degree 3 has chromatic number at least 8.

\begin{theorem}[$\chi_2(\mathcal{D}_3^c) \leq 10$]
	The class of connected $2$-edge-colored graphs with maximum degree $3$ has chromatic number at most $10$.
\end{theorem}

\begin{proof}
We proceed by contradicting the existence of a counter-example. Let $G$ be a 2-edge-colored graph such that $\chi_2(G) > 10$. \newline

\textit{Claim 1: $G$ contains no induced copy of $K_4^{s+}$ and no induced copy of $K_4^{s-}$.}

Assume otherwise. Let $SP_9^\dag$ be the 2-edge-colored graph formed from $SP_9$ by adding a new vertex $z$ so that there is a positive edge $zu$ for all $u\in \{0, 1, 2\}$ and a negative edge $zv$ for all $v \in \{{2x, 2x+1, 2x+2}\}$.

Let $v$ be a vertex connecting a $K_4^{s}$ (chosen arbitrarily) to the rest of the graph. Let $G'$ be obtained from $G$ after removing every $K_4^{s}$. Graph $G'$ is 2-degenerate and by Lemma~\ref{lem:3_2deg}, there exists a homomorphism $\varphi : G' \rightarrow SP_9$. We now extend $\varphi$ into a homomorphism $\varphi': G \rightarrow SP_9^\dag$. 

Figure~\ref{fig:3c_stitch+} shows how to color a $K_4^{s+}$ with $SP_9^\dag$ such that the vertex connecting the $K_4^{s+}$ to the rest of the graph is colored in  $x, x+1, x+2, 2x$ or $2x+1$. Note that every vertex of $SP_9$ is a positive (respectively negative) neighbor of at least one of this four vertices. We can therefore always extend the homomorphism to a $K_4^{s+}$. 

\begin{figure}[H]
  \begin{center}
    \subfloat{
    	\scalebox{0.8}
		{
			\begin{tikzpicture}[thick][thick]
        		\def \radius {1cm}
        		\def \margin {8} 
        		\tikzstyle{vertex}=[circle,minimum width=1em]
        		
        		\node[draw, vertex] (1) at (0,0) {};
        		\node[draw, vertex] (2) at (-1,-1) {};
        		\node[draw, vertex] (3) at (-1,1) {};
        		\node[draw, vertex] (6) at (-3,1) {};
        		\node[draw, vertex] (7) at (-3,-1) {};
        		
        		\draw (1.north) node[above]{$x/2x$};
        		\draw (2.north) node[below, shift={(0,-0.3)}]{$1$};
        		\draw (3.north) node[above]{$0$};
        		\draw (6.north) node[above]{$2$};
        		\draw (7.north) node[below, shift={(0,-0.3)}]{$z$};
        		
        		\draw[dashed] (1) -- (2);
        		\draw (1) -- (3);
        		\draw  (2) -- (6);
        		\draw  (2) -- (7);
        		\draw (3) -- (6);
        		\draw (3) -- (7);
        		\draw (6) -- (7);
        		\end{tikzpicture}
			
		}
    } 
    \subfloat{
    	\scalebox{0.8}
		{
			\begin{tikzpicture}[thick][thick]
        		\def \radius {1cm}
        		\def \margin {8} 
        		\tikzstyle{vertex}=[circle,minimum width=1em]
        		
        		\node[draw, vertex] (1) at (0,0) {};
        		\node[draw, vertex] (2) at (-1,-1) {};
        		\node[draw, vertex] (3) at (-1,1) {};
        		\node[draw, vertex] (6) at (-3,1) {};
        		\node[draw, vertex] (7) at (-3,-1) {};
        		
        		\draw (1.north) node[above]{$x+1/2x+1$};
        		\draw (2.north) node[below, shift={(0,-0.3)}]{$0$};
        		\draw (3.north) node[above]{$1$};
        		\draw (6.north) node[above]{$2$};
        		\draw (7.north) node[below, shift={(0,-0.3)}]{$z$};
        		
        		\draw[dashed] (1) -- (2);
        		\draw (1) -- (3);
        		\draw  (2) -- (6);
        		\draw  (2) -- (7);
        		\draw (3) -- (6);
        		\draw (3) -- (7);
        		\draw (6) -- (7);
        		\end{tikzpicture}
			
		}
    }
    \subfloat{
    	\scalebox{0.8}
		{
			\begin{tikzpicture}[thick][thick]
        		\def \radius {1cm}
        		\def \margin {8} 
        		\tikzstyle{vertex}=[circle,minimum width=1em]
        		
        		\node[draw, vertex] (1) at (0,0) {};
        		\node[draw, vertex] (2) at (-1,-1) {};
        		\node[draw, vertex] (3) at (-1,1) {};
        		\node[draw, vertex] (6) at (-3,1) {};
        		\node[draw, vertex] (7) at (-3,-1) {};
        		
        		\draw (1.north) node[above]{$x+2$};
        		\draw (2.north) node[below, shift={(0,-0.3)}]{$0$};
        		\draw (3.north) node[above]{$2$};
        		\draw (6.north) node[above]{$1$};
        		\draw (7.north) node[below, shift={(0,-0.3)}]{$z$};
        		
        		\draw[dashed] (1) -- (2);
        		\draw (1) -- (3);
        		\draw  (2) -- (6);
        		\draw  (2) -- (7);
        		\draw (3) -- (6);
        		\draw (3) -- (7);
        		\draw (6) -- (7);
        		\end{tikzpicture}
			
		}
    } 
    \caption{How to color a $K_4^{s+}$ with $SP_9^\dag$ such that the vertex connecting it to the rest of the graph is colored in $x, x+1, x+2, 2x$ or $2x+1$.} 
    \label{fig:3c_stitch+}
  \end{center}
\end{figure}
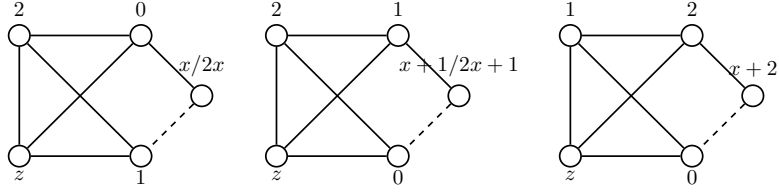

Similarly, Figure~\ref{fig:3c_stitch-} shows how to color a $K_4^{s-}$ with $SP_9^\dag$ such that  the vertex connecting the $K_4^{s-}$ to the rest of the graph is colored in $0, 1, x, x+1$ or $x+2$. Note that every vertex of $SP_9$ is a positive (respectively negative) neighbor of at least one of this four vertices. We can therefore always extend the homomorphism to a $K_4^{s-}$. 

\begin{figure}[H]
  \begin{center}
    \subfloat{
    	\scalebox{0.8}
		{
			\begin{tikzpicture}[thick][thick]
        		\def \radius {1cm}
        		\def \margin {8} 
        		\tikzstyle{vertex}=[circle,minimum width=1em]
        		
        		\node[draw, vertex] (1) at (0,0) {};
        		\node[draw, vertex] (2) at (-1,-1) {};
        		\node[draw, vertex] (3) at (-1,1) {};
        		\node[draw, vertex] (6) at (-3,1) {};
        		\node[draw, vertex] (7) at (-3,-1) {};
        		
        		\draw (1.north) node[above]{$0/x$};
        		\draw (2.north) node[below, shift={(0,-0.3)}]{$2x$};
        		\draw (3.north) node[above]{$2x+1$};
        		\draw (6.north) node[above]{$2x+2$};
        		\draw (7.north) node[below, shift={(0,-0.3)}]{$z$};
        		
        		\draw (1) -- (2);
        		\draw[dashed] (1) -- (3);
        		\draw[dashed]  (2) -- (6);
        		\draw[dashed]  (2) -- (7);
        		\draw[dashed] (3) -- (6);
        		\draw[dashed] (3) -- (7);
        		\draw[dashed] (6) -- (7);
        		\end{tikzpicture}
			
		}
    } 
    \subfloat{
    	\scalebox{0.8}
		{
			\begin{tikzpicture}[thick][thick]
        		\def \radius {1cm}
        		\def \margin {8} 
        		\tikzstyle{vertex}=[circle,minimum width=1em]
        		
        		\node[draw, vertex] (1) at (0,0) {};
        		\node[draw, vertex] (2) at (-1,-1) {};
        		\node[draw, vertex] (3) at (-1,1) {};
        		\node[draw, vertex] (6) at (-3,1) {};
        		\node[draw, vertex] (7) at (-3,-1) {};
        		
        		\draw (1.north) node[above]{$1/x+1$};
        		\draw (2.north) node[below, shift={(0,-0.3)}]{$2x+1$};
        		\draw (3.north) node[above]{$2x$};
        		\draw (6.north) node[above]{$2x+2$};
        		\draw (7.north) node[below, shift={(0,-0.3)}]{$z$};
        		
        		\draw (1) -- (2);
        		\draw[dashed] (1) -- (3);
        		\draw[dashed]  (2) -- (6);
        		\draw[dashed]  (2) -- (7);
        		\draw[dashed] (3) -- (6);
        		\draw[dashed] (3) -- (7);
        		\draw[dashed] (6) -- (7);
        		\end{tikzpicture}
			
		}
    }
    \subfloat{
    	\scalebox{0.8}
		{
			\begin{tikzpicture}[thick][thick]
        		\def \radius {1cm}
        		\def \margin {8} 
        		\tikzstyle{vertex}=[circle,minimum width=1em]
        		
        		\node[draw, vertex] (1) at (0,0) {};
        		\node[draw, vertex] (2) at (-1,-1) {};
        		\node[draw, vertex] (3) at (-1,1) {};
        		\node[draw, vertex] (6) at (-3,1) {};
        		\node[draw, vertex] (7) at (-3,-1) {};
        		
        		\draw (1.north) node[above]{$x+2$};
        		\draw (2.north) node[below, shift={(0,-0.3)}]{$2x+2$};
        		\draw (3.north) node[above]{$2x$};
        		\draw (6.north) node[above]{$2x+1$};
        		\draw (7.north) node[below, shift={(0,-0.3)}]{$z$};
        		
        		\draw (1) -- (2);
        		\draw[dashed] (1) -- (3);
        		\draw[dashed]  (2) -- (6);
        		\draw[dashed]  (2) -- (7);
        		\draw[dashed] (3) -- (6);
        		\draw[dashed] (3) -- (7);
        		\draw[dashed] (6) -- (7);
        		\end{tikzpicture}
			
		}
    } 
    \caption{How to color a $K_4^{s-}$ with $SP_9^\dag$ such that the vertex connecting it to the rest of the graph is colored in $0, 1, x, x+1$ or $x+2$.} 
    \label{fig:3c_stitch-}
  \end{center}
\end{figure}
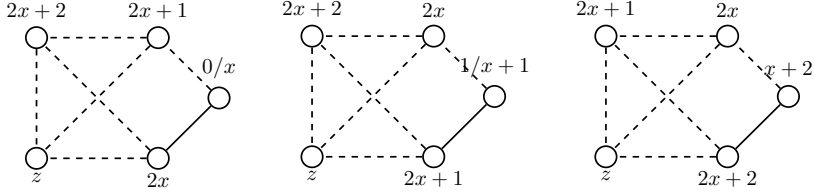

We can always find a 10-coloring of $G$, a contradiction.

\hfill$\Diamond$\newline

By Claim 1, $G$ cannot contain a $K_4^{s}$ hence we may assume that $G$ is 3-regular since if $G$ was not 3-regular we could find a 10-coloring by Lemma~\ref{lem:3_2deg}. \newline

\textit{Claim 2: $G$ contains no bridge.}

Assume otherwise. Let $uv$ be a bridge of $G$. Since $G$ cannot contain a copy of $K_4^s$, by Lemma~\ref{lem:3_2deg} there is a homomorphism $\varphi: G-uv \rightarrow SP_9$. Coloring $v$ in a 10th color yields a $10$-coloring of $G$, a contradiction.

\hfill$\Diamond$\newline

\textit{Claim 3: No vertex of $G$ is incident to three positive or three negative edges.}

Assume otherwise. Without loss of generality, let $v$ be a vertex of $G$ with neighbors $u_1, u_2, u_3$ so that each of $vu_1, vu_2, vu_3$ is positive. By Claim 1, G contains no $K_4^s$. Therefore, by Lemma~\ref{lem:3_2deg}, there is a homomorphism $\varphi: G - v \rightarrow SP_9$. We extend $\varphi$ to be a 10-coloring of $G$ by coloring $v$ in a 10th color.

\hfill$\Diamond$\newline

\textit{Claim 4: $G$ contains no copy of $K_3$.}

Assume otherwise. Let $u, v, w \in V(G)$ induce a copy of $K_3$ in $G$. Let $u'$ (resp. $v'$, $w'$) be the remaining neighbor of $u$ (reps. $v$, $w$). Without loss of generality, let $uv$ be negative (we can do this because $SP_9$ is antiautomorphic by Lemma~\ref{lem:PSP}). By Claim 1, $G$ contains no $K_4^s$. Form $G'$ from $G$ by removing the edge $uv$ and adding a vertex $z$, with positive edge $zu$ and negative edge $zv$. By Lemma~\ref{lem:3_2deg} there is a homomorphism $\varphi : G' \rightarrow SP_9$. Without loss of generality we may assume $\varphi(u) = 0$ and $\varphi(v) = 1$ by edge-transitivity of $SP_9$. If $\varphi(u') \neq 1$, coloring $u$ in a 10th color yields a 10-coloring of $G$. If $\varphi(v') \neq 0$, coloring $v$ in a 10th color yields a 10-coloring of $G$. We may now assume that $\varphi(u') = 1$ and $\varphi(v') = 0$. Consider now restricting $\varphi$ to $G - \{u, v, w\}$. We can extend $\varphi$ to a 10-coloring of $G$ by coloring $w$ in a 10th color and by letting $\varphi(u) = x+1$ and $\varphi(v) = 2x$ or $\varphi(u) = 2$ and $\varphi(v) = x$ such that $\varphi(u) \neq \varphi(w')$ and $\varphi(v) \neq \varphi(w')$, a contradiction. 

\hfill$\Diamond$\newline

By Claim~3, we may partition the vertices of $G$ in two sets $P$ and $N = V(G) \setminus P$ where vertices in $P$ are incident with exactly two positive edges and vertices in $N$ are incident with exactly two negative edges.\newline

\textit{Claim 5: There is no edge between a vertex of $P$ and a vertex of $N$.}

Assume otherwise. Consider $u \in P$ and $v \in N$. Without loss of generality let $uv$ be a negative edge. Let $u_1 \neq v$ and $u_2 \neq v$ be distinct neighbors of u. Since $u \in P$, the edges $uu_1$ and $uu_2$ are both positive. Note that by Claim~4, $G$ does not contain a copy of $K_3$ so $u_1$ and $u_2$ are not adjacent. Let $w$ be a neighbor of $v$ such that $vw$ is positive. Form $G'$ from $G$ by removing $u$ and adding a negative edge between $u_1$ and $u_2$. 

Note that a $K_4^s$ contains three vertices incident with only positive or negative edges. By Claim~3, $G$ does not contain such vertices. Adding the edge $u_1u_2$ to form $G'$ may create at most two vertices in $G'$ incident with 3 negative edges. Therefore, $G'$ does not contain a $K_4^s$.

Since $G$ does not contain a bridge by Claim~2, $G'$ is 2-degenerate and by Lemma~\ref{lem:3_2deg}, there is a homomorphism $\varphi : G' \rightarrow SP_9$. By Property $P_{2, 2}^*$ of $SP_9$ (Lemma~\ref{P22star}) we can extend $\varphi$ to include $u$ so that $\varphi(u) \neq \varphi(w)$. Coloring $v$ in a 10th color yields a 10-coloring of $G$, a contradiction.

\hfill$\Diamond$\newline

By Claim 5, we may assume that either $P$ or $N$ is empty, that is to say that either all the vertices are incident with exactly two positive edges or all the vertices are incident with exactly two negative edges.\newline
	
\textit{Claim 6: $G$ does not exist.} 

Without loss of generality we can restrict ourselves to the case in which every vertex is adjacent to exactly two positive edges. Let $u$ and $v$ be two vertices adjacent with a positive edge.
Let $G'$ be $G$ after removing edge $uv$ and adding a vertex $z$, with positive edge $zu$ and negative edge $zv$.
Graph $G'$ is 2-degenerate and by Claim~4, $G'$ contains no $K_3$ and therefore no $K_4^{s}$. By Lemma~\ref{lem:3_2deg}, there is a homomorphism $\varphi : G' \rightarrow SP_9$.
Let $u'$ be the negative neighbor of $u$. If $\varphi(u)\varphi(v)$ is a positive edge in $SP_9$ then $\varphi$ is already a homomorphism from $G$ to $SP_9$.
We may now assume that $\varphi(u)\varphi(v)$ is a negative edge in $SP_9$ since $u$ and $v$ cannot have the same color thanks to $z$. By the edge-transitivity of $SP_9$ we may assume without loss of generality that $\varphi(u) = 0$ and $\varphi(v) = x+1$.
Suppose that $\varphi(u') \neq x+1$ then extending $\varphi$ by coloring $u$ in a 10th color yields a 10-coloring of $G$.
We may now assume that $\varphi(u') = x+1$. By Property $P_{2, 2}^*$ of $SP_9$ (Lemma~\ref{P22star}), we can recolor $v$ with a different color which gives us $\varphi(u') \neq \varphi(v)$ and allows us to extend $\varphi$ into a 10-coloring of $G$ by coloring $u$ in a 10th color.
We can always find a 10-coloring of $G$, a contradiction.

\hfill$\Diamond$\newline
\end{proof}

\subsection{Signed graphs with maximum degree $3$}
\label{sec:3s}

Bensmail et al.~\cite{BPS} proved that every connected 2-edge-colored graph with maximum degree~3 except the all positive and all negative $K_4$ admits a homomorphism to $TR(SP_5)$, hence $\chi_2(\mathcal{D}_3^c)\le 12$, and $\chi_s(\mathcal{D}_3^c)\le 6$ by Lemma~\ref{lem:BG}. Their proof uses a computer to show that a minimal counter-example cannot contain some configurations and then concludes by using the properties of $TR(SP_5)$. 

In the non-connected case, we can easily get $\chi_2(\mathcal{D}_3)\le 14$ and thus $\chi_s(\mathcal{D}_3)\le 7$ by Lemma~\ref{lem:BG} (it is possible to create an all positive  $K_4$ and an all negative $K_4$ in $TR(SP_5)$ by adding two antitwinned vertices). 

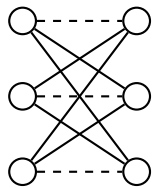
\begin{figure}[H]
	\centering
	\scalebox{1}
	{
		\begin{tikzpicture}[thick]
		\def \radius {1cm}
		\def \margin {8} 
		\tikzstyle{vertex}=[circle,minimum width=1em]
		
		\node[draw, vertex] (1) at (0, 0) {};
		\node[draw, vertex] (2) at (0, 1) {};
		\node[draw, vertex] (3) at (0, 2) {};
		\node[draw, vertex] (4) at (1.5, 0) {};
		\node[draw, vertex] (5) at (1.5, 1) {};
		\node[draw, vertex] (6) at (1.5, 2) {};
		
		\draw[dashed] (1) -- (4);
		\draw[dashed] (2) -- (5);
		\draw[dashed] (3) -- (6);
		
		\draw (1) -- (5);
		\draw (1) -- (6);
		\draw (2) -- (4);
		\draw (2) -- (6);
		\draw (3) -- (4);
		\draw (3) -- (5);
		\end{tikzpicture}
	}
	\caption{A signed clique with maximum degree 3 on 6 vertices.}
	\label{fig:3s_chi6}
\end{figure}

The following signed clique from Figure~\ref{fig:3s_chi6} with maximum degree 3 on 6 vertices gives us a lower bound of 6.

\bibliographystyle{unsrt}
\bibliography{biblio}

\end{document}